\documentclass[13pt,letterpaper,reqno]{article}
\usepackage{epsf}
\usepackage{amsmath}
\usepackage{amsfonts}
\usepackage{amssymb}
\usepackage{amstext}
\usepackage{amsbsy}
\usepackage{amsthm}
\usepackage{amscd}
\usepackage[all]{xy}

\newtheorem{theorem}{Theorem}[section]

\newtheorem{lem}[theorem]{Lemma}
\newtheorem{prop}[theorem]{Proposition}
\newtheorem{dfn}[theorem]{Definition}
\newtheorem{rem}[theorem]{Remark}
\newtheorem{cor}[theorem]{Corollary}
\newtheorem{claim}[theorem]{Claim}
\newcommand{\w}{\widetilde}

\newcommand{\bpf}{\noindent {\em Proof.  }}
\newcommand{\epf}{\qed \vspace{+10pt}}

\newcommand{\Zee}{\mathbb{Z}}
\newcommand{\Ree}{\mathbb{R}}
\newcommand{\Cee}{\mathbb{C}}

\newcommand{\eeight}{\left ( E_8 \times E_8 \right ) \rtimes \Zee_2}
\newcommand{\spint}{{\rm Spin}(32)/\Zee_2}
\newcommand{\mkt}{\mathcal{M}_{{\rm K3}}}
\newcommand{\lcsd}{\mathcal{V}}
\begin{document}
%
%  Title Stuff
%\pagestyle{plain}
%\title{
%{\large {On ${\bf K3}$ Surfaces with Large Complex Structure}}
%}
%\author[A. Clingher]{Adrian Clingher}
%\author[C. Doran]{Charles Doran}
%\address{Department of Mathematics, Stanford University, Stanford, CA 94305}
%\email{clingher@math.stanford.edu}
%\address{Department of Mathematics, University of Washington, Seattle, WA 98195}
%\email{doran@math.washington.edu}
%\date{Friday, June 18, 2004.}

\title{On K3 Surfaces with Large Complex Structure} 
%\footnote{Preliminary version. Temporary title.}}
\author{Adrian Clingher
\thanks{
Department of Mathematics, Stanford University, Stanford, CA 94305. {\bf e-mail:} {\it clingher@math.stanford.edu}
}
\and 
Charles F. Doran\thanks{
Department of Mathematics, University of Washington, Seattle, WA 98195. {\bf e-mail:} {\it doran@math.washington.edu}
}
} 
\maketitle
\begin{center}
\abstract{\noindent We discuss a notion of large complex structure 
for elliptic K3 surfaces with section inspired by the eight-dimensional 
F-theory/heterotic duality in string theory. This concept is naturally associated with 
the Type II Mumford partial compactification of the moduli space of 
periods for these structures. The paper provides an explicit Hodge-theoretic 
condition for the complex structure of an elliptic K3 surface with section to 
be large. We also establish certain geometric consequences of this 
large complex structure condition in terms of the Kodaira types of 
the singular fibers of the elliptic fibration.     }  
\end{center}
%
%
%
%
%
%\tableofcontents
%
\section{Introduction}
\noindent An elliptic K3 surface with section is a triple $(X,\varphi,S)$ consisting of a K3 surface $X$, 
an elliptic fibration $ \varphi \colon X \rightarrow \mathbb{P}^1 $ and a smooth rational curve $S$ making 
a section of $\varphi$. The extra structure given by the elliptic fibration and section on $X$ is equivalent 
to a pseudo-ample hyperbolic lattice polarization in the sense of Dolgachev \cite{dolgachev} and one can 
use this property to construct a coarse moduli space for elliptic K3 surfaces with section.
This space, which we shall denote here by $\mkt$, is a quasi-projective analytic variety of complex 
dimension eighteen. The properties of $\mkt$ have been extensively studied by means of the {\bf period map}:
\begin{equation}
\label{permap}
{\rm per} \colon \mathcal{M}_{{\rm K3}} \ \rightarrow \ \Gamma \backslash \Omega. 
\end{equation} 
which associates to a given triple $(X,\varphi,S)$ the equivalence class of its polarized Hodge structure. Here 
$ \Omega $ represents the classical domain of periods:
$$ 
\{ \ 
[\omega] \in \mathbb{P} \left( L \otimes \mathbb{C} \right ) \ \vert \ 
\omega \cdot \omega = 0, \ \omega \cdot \bar{\omega}  > 0  \ \} 
$$
where $L$ is the unimodular even lattice of signature $(2,18)$ and $\Gamma$ is the group 
of integral isometries of $L$ acting on $\Omega $ in a natural way. A special case of 
the Global Torelli theorem for lattice polarized $K3$ surfaces (see \cite{dolgachev}) asserts 
that $(\ref{permap})$ is an isomorphism of analytic spaces.   
\par The target space of the period map $(\ref{permap})$ is connected but not compact. 
However, due to the nice arithmetic features of the period domain, there exists
quite an array of methods at one's disposal for (partially) compactifying 
$\Gamma \backslash \Omega$. The simplest and most standard procedure is the Baily-Borel method
\cite{borel} which exploits the natural holomorphic identification between $\Omega$ and the hermitian 
symmetric space:
$$ {\rm O}(2,18)/{\rm SO}(2) \times {\rm O}(18), $$ 
in order to fully compactify $\Gamma \backslash \Omega$ 
by adding a number of curves and points. However, 
the Baily-Borel construction does not capture the full geometric information encoded in the periods, a fact 
reflected in the high codimension of the Baily-Borel boundary components. We shall be concerned here with a 
different compactification method, which realizes, to a certain extent, a blow-up of Baily-Borel's 
construction. This method is a special case of Mumford's
toroidal compactification construction \cite{ash}. The procedure, which was first applied in the 
K3 surface context by Friedman \cite{friedman1} \cite{frthesis}, constructs a smooth partial compactification 
\begin{equation}
\label{arithcomp}
\Gamma \backslash \Omega \  \subset \ \overline{\Gamma \backslash \Omega},
\end{equation}
by, essentially, adding to the quotient space $\Gamma \backslash \Omega$ two 
{\bf Mumford boundary divisors} $ \mathcal{D}_1 $ and $ \mathcal{D}_2 $ associated to 
the two distinct classes of Type II maximal rational parabolic 
subgroups of ${\rm O}(2,18)$. 
\par The Type II partial compactification $(\ref{arithcomp})$, although purely 
arithmetic in nature, has an interesting geometric interpretation on the moduli space side. 
It corresponds to an enlargement 
$$ \mkt \ \subset \ \overline{\mathcal{M}}_{{\rm K3}} $$
obtained by allowing certain normal crossing degenerations of elliptic $K3$ surfaces with section, 
the Type II stable elliptic $K3$ surfaces with section (see section 3 of \cite{clingher1}). 
The period map $ (\ref{permap})$ extends to an 
identification:
\begin{equation}
\overline{{\rm per}} \colon \overline{\mathcal{M}}_{{\rm K3}} \ \rightarrow \ 
\overline{\Gamma \backslash \Omega }
\end{equation}
and allows one to regard the boundary points of $ \mathcal{D}_1 $ and $ \mathcal{D}_2 $ 
as periods for the singular surfaces. 
\par One of the essential ingredients of this compactification method is the existence of 
a {\bf large complex structure domain} associated to each of the two boundary divisors. Let us postpone 
to section $\ref{rev}$ a more detailed definition of these domains and give just a brief description 
here. Consider a Type II maximal parabolic subgroup of 
${\rm O}(2,18)$ defined over the rationals and, denote by ${\rm P}$ its intersection with $\Gamma$.  
One has then a holomorphic non-normal covering projection with infinitely many sheets:
$$ \pi \colon {\rm P} \backslash \Omega \ \rightarrow \ \Gamma \backslash \Omega. $$
Moreover, the total space of this map fibers holomorphically, as in the diagram below, over the 
appropriate Mumford boundary component $ \mathcal{D}$, 
with all fibers being isomorphic to complex open punctured discs. 
\begin{equation}
\label{diag1}
\xymatrix{
\Omega \ \ar [r] & \ {\rm P} \backslash \Omega \ \ar [d] _{\alpha} \ar [r] ^{\pi} & \ \Gamma \backslash \Omega \\
& \mathcal{D}
}
\end{equation}   
One can select then in $\Omega$ a subset with special properties. 
\begin{theorem}
\label{plema}
There exists an open subset $ {\mathcal V} \subset \Omega $ such that: ${\mathcal V}$ is left invariant by the 
action of ${\rm P}$, the image of ${\mathcal V}$ under the natural projection to $ {\rm P} \backslash \Omega $ 
intersects each fiber of $\alpha$ over an open neighborhood of the puncture and, the $\Gamma$ equivalence reduces
to ${\rm P}$-equivalence on $\mathcal{V}$.
\end{theorem}
\noindent In the context of the above theorem, the restriction:
\begin{equation}
\label{gluemap}
\pi_{\vert {\rm P} \backslash \mathcal{V}} \colon {\rm P} \backslash \mathcal{V} \ 
\rightarrow \ \Gamma \backslash \Omega 
\end{equation}
is an isomorphism onto its range
\begin{equation}
\label{lcsd}
\mathcal{U} : = \pi \left ( {\rm P} \backslash \mathcal{V} \right ) .  
\end{equation} 
The inverse map of $(\ref{gluemap})$ provides then the 
essential gluing map which allows one to smoothly fit the boundary component $ \mathcal{D} $ together 
with  $\Gamma \backslash \Omega $. Moreover, one obtains a natural holomorphic fibration of 
$\mathcal{U}$ over the boundary divisor $\mathcal{D}$, whose fibers are copies of $\Cee^*$. 
The open subset $\mathcal{U} $ is the {\bf large complex structure domain} associated to $\mathcal{D}$. 
The inverse image of $ \mathcal{U}$ under the period map determines an open region of $\mkt$ which 
is said to correspond to elliptic K3 surfaces with {\bf large complex structure}.
\par The statement in Theorem $ \ref{plema} $ represents a special case of a general reduction 
theorem of Ash-Mumford-Rapoport-Tai proved in chapter 5 of \cite{ash}. However, the method of proof 
in \cite{ash} does not yield an explicit description of the open subset $\mathcal{U}$. It is therefore not an easy 
task to decide whether a given period line $[\omega]$ corresponds to an elliptic K3 surface with section of large complex structure or not. 
The goal of this paper is to introduce a simple, 
easy-to-test condition on the period lines in $\Omega$ leading to an 
open subset $ \mathcal{V} \subset \Omega $ which satisfies all requirements of Theorem 
$ \ref{plema}$. This, in addition to giving an alternative proof for Theorem $\ref{plema}$, 
provides an effective Hodge-theoretic method of testing whether a given elliptic K3 surface 
with section has large complex structure. 
\par Let us close this introductory section by also mentioning the string theory motivation 
underlying this work. This shall also serve as an explanation for the ``large complex structure''
terminology used for the subset $\mathcal{U}$ of $(\ref{lcsd})$. 
\par Following the works 
of Vafa \cite{vafa96} and Sen \cite{sen96} in 1996, 
it was noted that the geometry underlying elliptic $K3$ surfaces with 
section is related to the geometry of elliptic curves endowed with certain 
flat principal $G$-bundles. This non-trivial connection appears in string theory 
as the eight-dimensional 
manifestation of the phenomenon called F-theory/heterotic string duality. Over the past ten years 
the correspondence has been analyzed extensively (\cite{clingher1} \cite{curio} \cite{donagi} \cite{donagi2} \cite{fmw} \cite{vafa1} \cite{vafa2}) from a purely mathematical point of view. As it turns out, it leads to a beautiful 
geometric picture which links together moduli spaces for 
these two seemingly distinct types of geometrical objects: elliptic 
$K3$ surfaces with section and flat bundles over elliptic curves. 
\par In a brief description, what happens 
is the following. Let $G$ be one 
of the following two Lie groups: 
\begin{equation}
\label{picture}
\eeight \ \ \ \ \ \ \spint.
\end{equation}
As argued in \cite{friedmanmorgan}, one can define a moduli space $\mathcal{M}_{E,G}$ of equivalence classes of 
flat $G$-bundles over elliptic curves as a quasi-projective analytic space of complex dimension 
seventeen. There exists then a holomorphic isomorphism between  
$\mathcal{M}_{E,G}$ and one of the boundary divisors $ \mathcal{D}_i$ introduced in the previous 
Type II partial compactification of $ \Gamma \backslash \Omega $.  
Moreover, under this correspondence, the total space $ {\rm P} \backslash \Omega $ 
of the fibration $ \alpha $ in diagram $ (\ref{diag1})$ can be holomorphically identified with 
the moduli space $ \mathcal{M}_{{\rm het}} $ of 
classical vacua in heterotic string theory\footnote{$ \mathcal{M}_{{\rm het}} $ represents 
the moduli space of equivalence classes of pairs of flat $G$-bundles and complexified K\"{a}hler 
classes over elliptic curves. We refer the interested reader to \cite{clingher1} and \cite{clingher2} for 
more details regarding this space.}. 
\par The prediction made by the string duality is then that, although the two spaces $ \mkt $ and 
$ \mathcal{M}_{{\rm het}} $ are not globally identical, there should exist open regions in 
each of them, neighboring boundary divisors at infinity (regions that correspond, in physics language, 
to large levels of energy), that are analytically isomorphic. On the heterotic side, high energy levels 
appear when the two-torus has large volume and therefore, such a region 
has to be a tubular neighborhood of the punctures in the fibration $ \alpha $ of diagram $ (\ref{diag1})$. 
On the F-theory side, the appropriate region is not a priori obvious and, by convention, 
is said to correspond to elliptic K3 surfaces with section of large complex structure.   
\par Comparing the above paragraph with the arguments leading to the construction of the open subset 
$ \mathcal{U} $ in $ (\ref{lcsd}) $, we note that the statement of Theorem $ \ref{plema} $ precisely 
captures the feature predicted by the duality. Hence, one is naturally led to characterize the complex 
structures associated to elliptic K3 surfaces with section in the open region 
$$  \mathcal{U} \ \subset \ \Gamma \backslash \Omega $$ 
as large\footnote{This notion of large complex structure differs from the similarly-named condition arising in the 
context of Type IIA/IIB string duality (mirror symmetry) for K3 surfaces \cite{gross} \cite{morrison1}.}. 
\par The paper is organized as follows. In sections $\ref{prelim}$ and $ \ref{rev}$, we collect 
the basic facts needed to construct the period map and the Type II partial compactification 
of $\Gamma \backslash \Omega $. We introduce our Hodge-theoretic large complex structure condition 
in section $\ref{largecond}$. In section $\ref{proof} $, we show that the domain $\mathcal{V}$ 
defined by this 
condition satisfies the properties required by Theorem $ \ref{plema}$. Finally, in section 
$ \ref{ade} $, we discuss some geometric consequences of the large complex structure condition, 
in terms of what types of singular fibers can appear in the elliptic fibration associated to 
a large structure K3 surface.
\par The authors would like to thank John Morgan for helpful discussions throughout the evolution 
of this project.    
\section{Preliminaries}
\label{prelim}
A K3 surface $X$ is a non-singular, simply-connected complex surface with trivial canonical 
bundle. It is a well-known fact (see, for example, \cite{bpv}) that any two 
surfaces with these properties are diffeomorphic. 
The cohomology group $H^2(X,\Zee)$ is torsion free of rank $22$ and, when endowed with the pairing 
$ \langle \cdot, \cdot \rangle $ induced by the intersection form, it becomes an even unimodular 
lattice of signature $(3,19)$. There exists a unique lattice with these features. This 
lattice can be constructed independent of geometry by taking the orthogonal 
direct sum of the following irreducible factors:
\begin{equation}
\label{k3lat}
H \oplus H \oplus H \oplus E_8 \oplus E_8 
\end{equation}
where $H$ represents the rank-two hyperbolic lattice and $E_8$ is the unique negative-definite 
even and unimodular lattice of rank eight.            
\par The cohomology classes dual to algebraic cycles of $X$ span a special sublattice ${\rm NS}(X)$ of 
$H^2(X, \Zee)$, called the {\bf Neron-Severi lattice}. As a group, ${\rm NS}(X)$ is isomorphic to the 
Picard group of $X$, that is the group of algebraic equivalence classes of holomorphic line bundles over $X$. 
The rank the Neron-Severi lattice, denoted by $p_X$ varies between $0$ and $20$. By the Hodge index theorem, the signature of ${\rm NS}(X)$ is $(1,p_X-1)$. A generic K3 surface has rank $p_X=0$ and hence, is not projective.     
\par Our objects of interest are {\bf elliptically fibered K3 surfaces with section}. 
These are triples $(X, \varphi, S)$ consisting of a K3 surface $X$, a proper analytic map 
$\varphi \colon X \mapsto \mathbb{P}^1$ whose general fibers are smooth elliptic curves, 
and a smooth rational curve $S$ on $X$ which makes a section in the elliptic fibration. 
Two elliptically fibered $K3$ surfaces with section $(X, \varphi, S)$ and $(X', \varphi', S')$ are said to be 
{\bf equivalent} if there exists an analytic isomorphism $\alpha \colon X \mapsto X'$ with $\alpha(S)=S'$, 
inducing commutativity in the following diagram.
\begin{equation}
\label{diag11}
\xymatrix{
X \ar [dr] _{\varphi} \ar [rr] ^{\alpha} & & X' \ar [dl] ^{\varphi'} \\
& \mathbb{P}^1 & \\
}
\end{equation}
Given a triple $(X,\varphi,S)$ as above, one has two special classes $f,s \in {\rm NS}(X) $ 
associated to the elliptic fiber and section. These classes are independent and span 
a sublattice of rank two:  
\begin{equation}
\label{hyperemb}
\mathcal{H}_{(\varphi,S)} \subset {\rm NS}(X). 
\end{equation} 
In particular, $ p_X \geq 2$. The intersection form on $\mathcal{H}_{(\varphi,S)}$ with respect to the basis $\{f,s \} $ is 
\begin{equation}
\label{inters}
\left ( 
\begin{array}{cc}
0 & 1 \\
1 & -2 
\end{array}
\right ) 
\end{equation}
and therefore  $\mathcal{H}_{(\varphi,s)}$ is isometric to the standard rank-two hyperbolic lattice H. One has then a splitting of the Neron-Severi 
lattice of $X$ as an orthogonal direct sum:
\begin{equation}
{\rm NS}(X) \ = \ \mathcal{H}_{(\varphi,S)} \ \oplus \ \mathcal{W}_X 
\end{equation}
where $ \mathcal{W}_X$ is negative-definite and of rank $p_X-2$.

\vspace{.1in}
\begin{prop}
The sublattice $ \mathcal{H}_{(\varphi,S)} \subset {\rm NS}(X) $ totally determines the elliptic 
fibration with section $(\varphi,S)$ on $X$. 
\end{prop} 
\bpf
Let us assume that there exists a second elliptic structure with section $(\varphi',S')$ on $X$ 
such that:
\begin{equation}
\label{equalit}
 \mathcal{H}_{(\varphi,S)} \ = \ \mathcal{H}_{(\varphi',S')}.
 \end{equation}
Since the lattice $(\ref{equalit})$ is isometric to H, it contains only two classes of self-intersection 
$-2$. These two classes are $s$ and $-s$. By the Riemann-Roch theorem, only $s$ is effective. It 
follows that the two sections $S$ and $S'$ represent the same
$-2$ class and, since they are both irreducible curves it follows that $S=S'$. Next, we note that 
there exist only two isotropic elements in $(\ref{equalit})$ which have intersection 1 with $s$.
These classes are $f$ and $-f-s$. The second one cannot be effective and therefore, it must be that 
the two elliptic pencils $ \varphi$ and $\varphi'$ represent the same isotropic class $f$. Since 
the generic element in each of them is a smooth irreducible curve, it follows that $\varphi = \varphi'$.    
\epf

\vspace{.1in}
\noindent It is not true that all embeddings of H in $ {\rm NS}(X)$ are induced by 
an elliptic structure with section $(\varphi,S)$ on $X$. However, this statement becomes true if 
one requires that the embedded lattice contains a pseudo-ample class. 

\vspace{.1in}
\begin{dfn}
A class $d \in {\rm NS}(X)$ is called {\bf pseudo-ample} if it represents an effective divisor on $D$ which 
is nef (has non-negative intersection with any effective class) and has positive self-intersection. 
\end{dfn}

\vspace{.1in}
\noindent The above terminology is closely related to the classical notion of ampleness. Given a pseudo-ample 
class on $X$ and $D$ an effective divisor representing $d$, a result of Mayer \cite{mayer} asserts that 
the linear system $\vert nD \vert $ is base point free 
for $n \geq 3 $ and the associated map:
$$ \psi_{\vert nD \vert } \colon X \rightarrow \mathbb{P}^N $$
is a birational morphism. Moreover, the image of $\psi_{\vert nD \vert }$ is the normal model of $X$ obtained by contracting all curves not met 
by $D$ which are rational double point configurations.

\vspace{.1in}
\begin{theorem}
\label{connectdolg}
Let $ \mathcal{H} \subset {\rm NS}(X)$ be a sublattice isometric to H. There exists an elliptic structure 
with section $(\varphi,S)$ on $X$ such that:
$$ \mathcal{H}_{(\varphi,S)} \ = \ \mathcal{H} $$ 
if and only if $ \mathcal{H} $ contains a pseudo-ample class.  
\end{theorem}
\bpf
Let us first check that the above condition is necessary. Given an elliptic structure with section 
$(\varphi,S)$, one has the special classes $f,s \in \mathcal{H}_{(\varphi,S)}$. The set of effective 
classes on $X$ is then the semi-group generated by $f$, the irreducible components of the singular 
fibers, and the possible (multi)-sections. It is then a simple verification to check that the class $d=2f+s$ is
numerically effective and has positive self-intersection. It follows that $d$ is pseudo-ample.
\par In order to prove the opposite implication, we employ the following strategy. Assume that 
$\mathcal{H}$ is a sublattice of ${\rm NS}(X)$ which is isometric to H and contains a 
pseudo-ample class $d$. Denote by 
$\Gamma_X$ the group of isometries of $H^2(X,\Zee)$ whose $\Cee$-linear extensions preserve the 
Hodge filtration of $X$. We shall 
construct an elliptic fibration with section $(\varphi,C_0)$ such that:
\begin{equation}
\label{latticetransf11}
\mathcal{H} \ = \ \beta \left ( \mathcal{H}_{(\varphi,C_0)} \right ) 
\end{equation}  
for some $\beta \in \Gamma_X $. Then, by using standard arguments from the Strong Torelli Theorem for K3 surfaces, we argue that $\beta$ is induced by an analytic 
automorphism of $X$.
\par First, note that we can assume that $d$ has self-intersection $2$. There exist then two distinct 
classes in $\mathcal{H}$ of self-intersection $-2$. They differ by a change in sign. By the Riemann-Roch theorem, 
one and only one of these two classes is effective. Denote this class by $s$. Next, let $f$ be the 
unique isotropic element 
of $\mathcal{H}$ such that $d=2f+s$. By a result of Pjatecki\u{i}-\v{S}apiro and \v{S}afarevi\v{c} 
(see chapter 3 of \cite{shapiro}), there exists 
an isometry $\beta_1$ of $H^2(X,\Zee)$ such that $\beta_1(f)$ is the class of a smooth elliptic curve $F$ inducing 
an elliptic pencil $$ \varphi \colon X \rightarrow \mathbb{P}^1. $$ Moreover, the isometry $\beta_1$ 
is a composition of reflections associated to effective $-2$ classes in ${\rm NS}(X)$, and therefore $\beta_2 \in \Gamma_X$.  
\par We construct a section for the elliptic pencil $ \varphi $. The Riemann-Roch 
theorem combined with the fact that 
$$ [F] \cdot \beta_1(s) = 1, $$ 
implies that $ \beta_1(s)$ is 
an effective class. Therefore $ \beta_1(s)$ can be associated to a formal sum of 
irreducible curves:
$$ C \ = \ \sum_{i} \ n_i C_i , \ \ n_i \geq 1. $$  
But $ C_i \cdot F \geq 0 $ and $C \cdot F = 1 $. Hence, among the irreducible curves $C_i$ there 
exists a unique one, say, denoted by $C_0$, such that its intersection pairing with $F$ is not zero. 
We see that
$F \cdot C_0 = 1$, $n_0=1$, and $F \cdot C_i =0$ for $i\neq 0$. The restriction 
of the elliptic pencil $\varphi$ to $C_0$ then defines a degree one map $ C_0 \rightarrow \mathbb{P}^1$ 
and, from this we conclude that $C_0$ is a smooth rational curve. We have therefore a 
elliptic fibration with section $(\varphi, C_0)$ on $X$. 
\par Let then $ \mathcal{H}' : = \ \beta_1 \left ( \mathcal{H} \right ) $. The lattice $\mathcal{H}'$ 
is spanned by $[F]$ and $[C]$ and induces an orthogonal direct sum decomposition:
$$ {\rm NS}(X) \ = \ \mathcal{H}' \ \oplus \ \mathcal{W}. $$
In this decomposition, one can write:
$$ [C_0] \ = \ - \frac{(w_o,w_o)}{2}[F] + [C] + w_o $$
for some fixed $ w_o \in \mathcal{W}. $ Define then:
$$ \beta_2 \left ( a [F] + b[C] + w \right ) \ = \  
\left ( a - b \frac{(w_o,w_o)}{2} - (w, w_o) \right ) [F] + b [C] + bw_o + w.      $$  
One can see that $\beta_2$ is an isometry of $ {\rm NS}(X) $ satisfying $ \beta_2([F])=[F]$ and 
$ \beta_2([C])=[C_0]$. Moreover, $\beta_2$ is the restriction 
of an isometry in $\Gamma_X$ which, by an abuse of notation, we shall also name $\beta_2$. 
\par Let then 
$$\beta = \left ( \beta_2 \circ \beta_1 \right )^{-1}. $$ 
Since $ \beta([C_0])=s $ and $ \beta ([F])=f$, we find that $(\ref{latticetransf11})$ holds. We recall then 
the following result of Pjatecki\u{i}-\v{S}apiro and \v{S}afarevi\v{c} (Theorem 1 in chapter 6 of \cite{shapiro}):
\begin{theorem} (Pjatecki\u{i}-\v{S}apiro, \v{S}afarevi\v{c})
One has a decomposition:
$$ \Gamma_X \ = \ \Gamma^{{\rm eff}}_X \cdot {\rm W}(X) \cdot \{ \pm {\rm id} \} $$
where $\Gamma^{{\rm eff}}_X$ is the subgroup of effective isometries\footnote{By definition, an isometry 
$\gamma \in \Gamma_X $ is effective if it preserves the set of effective classes of $X$.} of $\Gamma_X$ and 
${\rm W}(X) $ is the subgroup generated by reflections with respect to effective $-2$ classes in ${\rm NS}(X)$. 
\end{theorem}   
\noindent One can write then 
$$\beta = \beta_1 \circ \beta_2 \circ \beta_3 $$
with $\beta_1 \in \Gamma^{{\rm eff}}_X $, 
$ \beta_2 \in {\rm W}(X) $ and $ \beta_3 = \pm {\rm id} $. Since $\beta(2[F]+[C_0]) = 2f+s = d $, we deduce that 
$\beta $ preserves the positive cone of $X$. By definition, so do $\beta_1$ and $\beta_2$. It follows therefore 
that $\beta_3 = {\rm id}$. 
\par Denote by $\mathcal{C}^+_X$ the K\"{a}hler cone of $X$ and let $ \overline{\mathcal{C}}^+_X$ be its closure 
inside the positive cone. Since $d$ and $2[F]+[C_0]$ are pseudo-ample classes, they both belong to $ \overline{\mathcal{C}}^+_X$. Therefore, $\beta$ sends an element of  
$ \overline{\mathcal{C}}^+_X$ to another element of $ \overline{\mathcal{C}}^+_X$. 
\par Recall then the following standard facts (\cite{bpv}, \cite{shapiro}). Isometries in 
$\Gamma^{{\rm eff}}_X$ preserve the K\"{a}hler cone (Proposition 3.10 in chapter VIII of \cite{bpv}) and therefore 
they preserve $ \overline{\mathcal{C}}^+_X$. On the other hand, $ \overline{\mathcal{C}}^+_X$ is a fundamental 
domain for the action of ${\rm W}(X)$, in the sense that any ${\rm W}(X)$-orbit meets 
$ \overline{\mathcal{C}}^+_X$ in exactly one point (Proposition 3.9 in chapter VIII of \cite{bpv}).
\par The last two paragraphs imply then that $\beta_2= {\rm id}$ and therefore $\beta \in \Gamma^{{\rm eff}}_X$. 
But then, by the Strong Torelli Theorem for K3 surfaces (see Theorem 11.1 in chapter VIII of \cite{bpv} or the 
similar results in \cite{shapiro} and \cite{looijenga}), there exists an analytic automorphism $\alpha \in {\rm Aut}(X)$ such that $\beta = \alpha^*$. It follows then that $(\varphi \circ \alpha, \alpha^{-1}(C_0)) $ is an 
elliptic fibration with section on $X$ and:
$$ \mathcal{H}_{(\varphi \circ \alpha, \alpha^{-1}(C_0))} \ = \ 
\alpha^* \left ( \mathcal{H}_{(\varphi,C_0) }\right ) \ = \ \mathcal{H}. $$ 
\epf 

\vspace{.1in}
\noindent In \cite{dolgachev}, Dolgachev has introduced the notion of {\bf pseudo-ample lattice polarization} of a K3 surface $X$. Given an even lattice M of signature $(1,t)$, $t\leq p_X-1$, which can be embedded in the K3 lattice, a pseudo-ample 
{\bf M-polarization} of $X$ is a lattice 
embedding 
$$ i \colon {\rm M} \hookrightarrow {\rm NS}(X) $$ 
whose image contains a pseudo-ample class. In this context, if H is the standard rank-two hyperbolic lattice , Theorem $\ref{connectdolg}$ shows that there  
exists a bijective correspondence    
\begin{equation}
\label{diag223}
\xymatrix{
\left \{ \ 
\txt{elliptic fibrations\\ 
with section on $X$ } \
\right \} 
\ \ar @{<->}  [r] & \    
\left \{ \
\txt{pseudo-ample   \\ 
H-polarizations of $X$} 
\ \right \} 
}. 
\end{equation} 
Moreover, under the above correspondence the equivalence relation on elliptic fibrations with section defined in $(\ref{diag11})$ translates precisely to Dolgachev's notion of equivalence for lattice polarizations. This leads 
to a canonical bijective correspondence between the equivalence classes of the two structures. One can therefore 
obtain a moduli space for triples $(X,\varphi,S)$ by constructing a moduli space for pairs $(X,\mathcal{H})$ 
consisting of a K3 surface X and a pseudo-ample H-polarization $\mathcal{H}$.   
\par The construction of such a moduli space of lattice polarizations has been done in \cite{dolgachev}.
We shall present just the main features. 
\par Given a triple $(X,\varphi,S)$, the orthogonal complement:
\begin{equation}
\label{k33flat}
\left ( \mathcal{H}_{(\varphi,S)} \right ) ^{\perp} \ \subset \ H^2(X, \Zee) 
\end{equation}
is an even, unimodular lattice of signature $(2,18)$, and therefore, by standard lattice theory, 
it is isometric to:
\begin{equation}
\label{k3flat}
{\rm L} \ = \ {\rm H} \oplus {\rm H} \oplus {\rm E}_8 \oplus {\rm E}_8. 
\end{equation} 
An isometry between $ (\ref{k33flat})$ and $ (\ref{k3flat})$ is called a marking. 
If two elliptic $K3$ surfaces with section $(X, \varphi,S)$ and $(X', \varphi',S')$ are given markings 
$q$ and $q'$, an analytic isomorphism $ \alpha \colon X \rightarrow X'$ preserving the elliptic fibrations and sections is said to be compatible with the markings if the induced isometry $\alpha^*$ between  
$  \left ( \mathcal{H}_{(\varphi',S')} \right ) ^{\perp} $ and $ \left ( \mathcal{H}_{(\varphi,S)} \right ) ^{\perp} $
satisfies $ q' = q \circ (\alpha^*)^{\perp} $.
\par One defines then the {\bf H-polarized period domain}:
\begin{equation}
\Omega \ = \ 
\{ \ 
[\omega] \in \mathbb{P}^1 \left ( {\rm L} \otimes \Cee \right ) \ \vert \ 
\langle \omega, \omega \rangle  =0, \ \langle \omega, \bar{\omega} \rangle  > 0  
\ \}.
\end{equation}
Every equivalence class of a marked elliptic K3 surface with section $(X,\varphi,S,q)$ determines uniquely a {\bf period line} 
\begin{equation}
\label{percoresp}
[\omega] \ = \ q \left ( H^{0,2}(X) \right )  \in \Omega. 
\end{equation}
This period correspondence can be seen naturally as an analytic morphism. Indeed, there exists a 
fine moduli space $ \mathcal{M}^{{\rm marked}}_{K3} $ for marked elliptic K3 surfaces with section. 
$\mathcal{M}^{{\rm marked}}_{K3} $ is an analytic space of complex dimension $18$ and 
it is constructed, along the lines of \cite{asterix} and \cite{shapiro}, by gluing together local 
moduli spaces of marked elliptic K3 surfaces with section. A lattice polarized version of 
the Global Torelli Theorem (\cite{burns}, \cite{shapiro} \cite{todorov}) then asserts that 
the period correspondence:
\begin{equation}
\label{permapmm}
\mathcal{M}_{{\rm K3}}^{{\rm marked}} \ \rightarrow \ \Omega 
\end{equation}
defined by $(\ref{percoresp})$ is a surjective morphism of analytic spaces.
\par This picture can be further refined by removing the markings. The discrete 
group $\Gamma $ of integral isometries of ${\rm L}$ acts on $\Omega$ as well as on 
$ \mathcal{M}^{{\rm marked}}_{{\rm K3}}$ and the morphism $(\ref{permapmm})$ is equivariant 
with respect to the two actions. The quotient space:
$$ \mkt := \ \Gamma \backslash \mathcal{M}^{{\rm marked}}_{{\rm K3}} $$ 
is a coarse 
moduli space for K3 surfaces with pseudo-ample H-polarizations (see remark 3.4 of 
\cite{dolgachev}). Then, by taking into account Theorem $\ref{connectdolg}$, one can show 
that $\mkt$ can be seen as a coarse 
moduli space for elliptic K3 surfaces with section. The induced {\bf period map}:
$$ {\rm per} \colon \mkt \ \rightarrow \ \Gamma \backslash \Omega $$ 
is an isomorphism of analytic spaces.       
\par The quotient $\Gamma \backslash \Omega$ is a connected space. To see the this, 
note that $\Omega$ is an open subset of an $18$-dimensional complex quadric and consists of two 
connected components interchanged by complex conjugation. But, the discrete group $\Gamma$ acts 
transitively on the set of these two connected components and hence the quotient space is connected.     
\par The space $\Gamma \backslash \Omega$ is, however, not compact. There are nonetheless methods 
available to compactify this space (and hence the moduli space $\mkt$), due to the fact that $\Gamma \backslash \Omega$ can be 
identified with an Hermitian symmetric space factored by the action of an arithmetic group. 
Indeed, the real Lie group ${\rm O}(2,18)$ of real isometries of ${\rm L} \otimes \mathbb{R}$ acts 
transitively on $ \Omega $. The action leads to an identification:
\begin{equation}
\label{ident}
\Omega \ = \ {\rm O}(2,18)/{\rm SO}(2) \times {\rm O}(18).
\end{equation}
The discrete group $\Gamma$ has an obvious action on the right term above, 
and the correspondence $(\ref{ident})$ becomes 
$\Gamma$-invariant. Therefore the above identification 
can be pushed to the level of quotients, where one obtains:
\begin{equation}
\Gamma \backslash \Omega \ = \ \Gamma \backslash {\rm O}(2,18)/{\rm SO}(2) \times {\rm O}(18). 
\end{equation}
In this context, a special case of Mumford's toroidal compactification \cite{ash} allows one 
to enlarge $ \Gamma \backslash \Omega $ by adding boundary components associated to maximal 
rational parabolic subgroups of Type II in ${\rm O}(2,18)$. These groups are, 
essentially, subgroups of isometries stabilizing a given rank-two primitive and isotropic 
sublattice ${\rm V} \subset {\rm L}$. 

\section{Review of the Type II Partial Compactification}
\label{rev}
One performs the Type II partial compactification of $ \Gamma \backslash \Omega $ by 
adding two specific boundary divisors. Let us give here a brief account of the procedure. 
For a detailed presentation we refer the reader to section 3 of \cite{clingher1} as well as to
\cite{friedman1} (for a closely related case). 
\par As mentioned at the end of the previous section, the procedure we are about to describe 
is centered around rational Type II parabolic subgroups of ${\rm O}(2,18)$ which, in turn, are 
defined as stabilizer groups for rank-two primitive and isotropic sublattices ${\rm V} $ of ${\rm L}$. 
By classical lattice theory \cite{nikulin}, there are two distinct types 
of sublattices with these features. Let ${\rm I}_2({\rm L})$ be the set of all possible such 
sublattices. The group $\Gamma$ acts on ${\rm I}_2({\rm L})$ and the action generates two 
distinct orbits. One can differentiate between these two orbits 
by analyzing the isomorphism class of the induced quotient lattice ${\rm V}^{\perp}/{\rm V}$.
This lattice is always negative-definite, even, unimodular and of rank $16$. There are only two 
non-equivalent lattices with these features. One is ${\rm E}_8 \oplus {\rm E}_8 $, the orthogonal 
direct sum of two copies of the unique negative definite, unimodular and even lattice of rank $8$. 
The other is the Barnes-Wall lattice $ {\rm D}^+_{16}$. 
\par Let us then outline the compactification construction. As a first step, pick one of the 
connected components of $\Omega$ and denote it by $\Omega^+$. One has then a natural identification: 
\begin{equation}
\Gamma \backslash \Omega \ = \ \Gamma^+ \backslash \Omega^+ 
\end{equation}
where $\Gamma^+$ is the index-two subgroup of $ \Gamma$ corresponding to isometries preserving 
$\Omega^+$. Select then a sublattice ${\rm V}$ with the properties described above. As the next step, 
the technical ingredient required by Mumford's method is a primitive integral 
and nilpotent element ${\rm N}$ in the Lie algebra of ${\rm O}(2,18)$ such that 
$ {\rm Im}({\rm N}) = {\rm V}$. At this point in the general construction of \cite{ash}, one has to make a 
choice of such ${\rm N}$. However, in the present situation, the choice is canonical. An 
endomorphism ${\rm N}$ with above characteristics is unique, up to a sign change. Moreover,
given such an ${\rm N}$, for any $[\omega] \in \Omega $, the quantity  
$ i \langle {\rm N}\omega , \bar{\omega}\rangle $ is real and non-zero. This quantity is positive 
on one connected component of $\Omega$ and negative on the other. We can then canonically 
select ${\rm N}$ such 
that $ i \langle {\rm N}\omega , \bar{\omega}\rangle >0 $ for $[\omega]$ in $\Omega^+$.
\par We next introduce the following group:
$$ {\rm U}({\rm N})_{\Cee} \ = \ \{ \ {\rm exp}(z{\rm N}) \ \vert \ z \in \Cee \ \}.  $$
This group is by definition isomorphic to $ (\Cee, +) $ and acts naturally 
on the compact quadric:
$$
\Omega^{\vee} \ 
= \ 
\{ \ 
[\omega] \in \mathbb{P}^1 \left ( {\rm L} \otimes \Cee \right ) \ \vert \ 
\langle \omega, \omega \rangle  =0  
\ \}
$$
in which $ \Omega $ embeds as an open subset. Denote 
by ${\rm U}({\rm N})_{\Zee} $ its subgroup  corresponding to $ z \in \Zee $ and set:
$$\Omega^+(V) = {\rm U}({\rm N})_{\Cee} \cdot \Omega^+ . $$ 
This allows one to construct 
the following projection:
\begin{equation}
\label{proj22}
{\rm U}({\rm N})_{\Zee} \backslash \Omega^+(V) \ \rightarrow \ {\rm U}({\rm N})_{\Cee} 
\backslash \Omega^+(V).  
\end{equation}
It can be easily verified that $(\ref{proj22})$ is a holomorphic principal bundle with 
structure group $ \mathbb{T} = {\rm U}({\rm N})_{\Zee} \backslash {\rm U}({\rm N})_{\Cee} $ and 
that ${\rm U}({\rm N})_{\Zee} \backslash \Omega^+ $ embeds as an open subset of the 
total space ${\rm U}({\rm N})_{\Zee} \backslash \Omega^+(V)$.
\par Let then ${\rm P} \subset \Gamma $ be the parabolic subgroup associated to ${\rm V}$, that is, 
the subgroup of isometries stabilizing ${\rm V}$. Denote ${\rm P}^+ = {\rm P} \cap \Gamma^+$. One checks that 
$  {\rm U}({\rm N})_{\Zee} \vartriangleleft {\rm P}^+$. Moreover, the discrete group ${\rm P}^+$ acts on both the 
base space and total space of the principal bundle $(\ref{proj22})$ and the action is compatible with both 
the projection, and the action of the structure group $\mathbb{T}$ on the total space $
{\rm U}({\rm N})_{\Zee} \backslash \Omega^+(V)$. Therefore,  $(\ref{proj22})$ descends to a holomorphic 
principal bundle\footnote{
%
% FOOTNOTE
Strictly speaking, the fibration $(\ref{proj333})$ is not a principal bundle, as it is not locally trivial and the 
group $\mathbb{T}$ does not act freely on certain fibers. A more rigorous description 
here would be that $(\ref{proj333})$ is a holomorphic Seifert $\mathbb{T}$-fibration.
} 
of structure group $ \mathbb{T}$:
\begin{equation}
\label{proj333}
\alpha_V \colon {\rm P}^+ \backslash \Omega^+({\rm V}) \ \rightarrow \ {\rm P}^+ \backslash \left ( 
{\rm U}({\rm N})_{\Cee} \backslash \Omega^+({\rm V})
\right ). 
\end{equation} 
For simplicity, we shall denote the base space of $(\ref{proj333})$ by $\mathcal{D}(V)$. This is the 
{\bf Mumford boundary component associated to ${\rm V}$} that one makes use of in order to partially compactify 
$ \Gamma^+ \backslash \Omega^+$. 
\par Mumford's toroidal compactification idea applies here in the following form. Since the 
group $ \mathbb{T} $ is naturally identified 
with $ \Cee^* $ by the exponential map, one can make this group act on a copy of the complex plane $\Cee$ in the 
standard way. This action is then used to construct the associated fibration with lines:
\begin{equation}
\label{comp22}
\overline{
{\rm P}^+ \backslash \Omega^+(V) 
} \ : = \  
{\rm P}^+ \backslash \Omega^+(V)  \times _{\mathbb{T}} \Cee \ ,
\end{equation}
a procedure that has the effect of adding to the total space of $(\ref{proj333})$ a 
divisor corresponding to the compactification of $\Cee^*$ to $\Cee$ in each fiber. With this construction 
in place, one defines:
\begin{equation}
\overline{
{\rm P}^+ \backslash \Omega^+ 
} \ : = \ {\rm interior \ of \ the \ closure \ of}  \ {\rm P}^+ \backslash \Omega^+ \ {\rm in} \ \overline{
{\rm P}^+ \backslash \Omega^+(V) 
}.
\end{equation}
Set-theoretically,  
$$ 
\overline{
{\rm P}^+ \backslash \Omega^+ 
}  \ = \ {\rm P}^+ \backslash \Omega^+ \ \ \bigsqcup \ \ \mathcal{D}({\rm V})
$$
and therefore we have just performed a holomorphic gluing of the boundary component 
$\mathcal{D}({\rm V})$ to $ {\rm P}^+ \backslash \Omega^+$. 
\par The last step in the procedure is the attaching of $\mathcal{D}({\rm V})$ to $\Gamma^+ \backslash \Omega^+$. 
The key ingredient here, as mentioned in the introduction, is the existence of a large complex structure period 
domain. That is, there exists an open subset: 
\begin{equation}
\label{lcsreq}
\lcsd^+ \ \subset \  \Omega^+ 
\end{equation} 
that satisfies the following properties: 
\begin{itemize}
\item [(a)] $\lcsd^+ $ is invariant under the action of ${\rm P}^+$,
\item [(b)] the restriction of ${\rm U}({\rm N})_{\Zee} \backslash \lcsd^+ $ to any given fiber 
of $(\ref{proj22})$ is an open neighborhood of the cusp,
\item [(c)] on $\lcsd^+$, the equivalence under the action of $\Gamma^+$ reduces to ${\rm P}^+$-equivalence.
%given any two $ [\omega_1], [\omega_2] \in \mathcal{V} $ and $\gamma \in \Gamma^+$ such that 
%$ \gamma([\omega_1]) = [\omega_2] $, one has that $ \gamma \in {\rm P}^+ $.  
\end{itemize}
If these three conditions are satisfied then the non-normal covering projection  
$ \pi \colon {\rm P}^+ \backslash \Omega^+ \rightarrow \Gamma^+ \backslash \Omega^+ $, when restricted 
to ${\rm P}^+ \backslash \lcsd^+$ induces a holomorphic isomorphism between ${\rm P}^+ \backslash \lcsd^+$ 
and its image. The following commutative diagram:
\begin{equation}
\label{bigdiag}
\xymatrix{
{\rm P}^+ \backslash \lcsd^+ \ \ar @{^{(}->} [r] \ar [d] _{\simeq} & \ 
{\rm P}^+ \backslash \Omega^+ \ar [r] \ar [d] _{\pi} & 
\overline{{\rm P}^+ \backslash \Omega^+} \\
\pi \left ( {\rm P}^+ \backslash \lcsd^+ \right ) \ \ar @{^{(}->} [r] & \ \Gamma^+ \backslash \Omega^+\\
}
\end{equation}
allows one to use (the inverse of) this isomorphism as gluing map, thus attaching $\mathcal{D}(V)$ holomorphically onto $ \Gamma^+ \backslash \Omega^+$.
\par The open set 
$$ \mathcal{U} = \pi \left ( {\rm P}^+ \backslash \lcsd^+ \right ), $$ 
is the {\bf large complex structure domain} associated to $\mathcal{D}({\rm V})$. 
\begin{rem}
Note that the subset $ \mathcal{U} $, as well as the gluing map,  
depends only on the $\Gamma$-orbit of ${\rm V}$ in ${\rm I}_2({\rm L})$. If one repeats the 
above procedure starting with a different sublattice $V' = \gamma(V) $ for $\gamma \in \Gamma $
the construction produces the same partial compactification. However, isotropic sublattices $V$ and 
$V'$ inequivalent under $\Gamma$ lead to distinct large complex structure regions and different 
compactifications.
\end{rem}
\noindent Finally, the Type II partial compactification $ \overline{\Gamma \backslash \Omega } $ is the space 
obtained after gluing onto $ \Gamma \backslash \Omega $ the two boundary components $\mathcal{D}_1$ 
and $ \mathcal{D}_2 $ associated to the two $\Gamma$-orbits of ${\rm I}_2({\rm V})$. It follows 
that $ \overline{\Gamma \backslash \Omega } $ is a quasi-projective analytic space of complex 
dimension eighteen. It contains $ \Gamma \backslash \Omega $ as a Zariski open subset. The complement
$$ \overline{\Gamma \backslash \Omega } \ \setminus \ \Gamma \backslash \Omega $$
is the disjoint union of two irreducible divisors. Each of the two divisors 
$\mathcal{D}_1$, $\mathcal{D}_2$ is a quotient of a smooth space by a finite group action. We refer 
the reader to \cite{ash} for proofs of these statements. 

\section{A Large Complex Structure Condition} 
\label{largecond}
The main goal of this paper is to introduce an effective Hodge-theoretic condition for an elliptic 
$K3$ surface with section to have large complex structure. In other words, we wish to construct 
explicitly a large complex structure period region $ \lcsd^+$ as in $(\ref{lcsreq})$ leading to 
a large complex structure domain $ \mathcal{U} = \pi \left ( {\rm P}^+ \backslash \mathcal{V} \right )$. 
%In this section, we shall present the condition.
%The proof that this condition complies with the large complex structure requirements discussed 
%in section $\ref{rev}$ will be given in section TBA. 
\par Let $(X, \varphi,S)$ be an elliptically fibered $K3$ surface with section. Assume that $V$ is a 
rank-two primitive and isotropic sublattice of 
$$ V \subset \ \left ( \mathcal{H}_X \right )^{\perp} \ \subset \ H^2(X, \Zee). $$ 
We shall attach to such a quadruple, $(X,\varphi,S,V)$, two Hodge-theoretic 
quantities denoted (by abuse of notation) $ \tau(X,V)$ and $ \widetilde{u}_2(X,V) $. The first quantity 
is a complex number with positive imaginary part. The second one is a positive real number.  
%The meaning of these two parameters will be clear 
%in section $\ref{narainsect}$ where they will appear as coordinates on the period domain.
\par Consider $ \omega \in H^2(X, \Cee) $ to be a class representing a holomorphic two-form on $X$. 
The class $\omega$ is unique up to a scaling by a non-zero complex number. Moreover, since the 
lattice $ \mathcal{H}_X$ is spanned by classes representing algebraic cycles, $ \omega $ belongs to 
\begin{equation}
\label{latperp}
\left ( \mathcal{H}_X \right )^{\perp} \otimes  \Cee.  
\end{equation} 
The intersection form on the real version of $ (\ref{latperp})$ has 
signature $(2,18)$, and hence (see also Lemma $\ref{lem21}$), $\langle \omega , y \rangle \neq 0 $ 
for any $y \in V$. We select then an oriented basis $\{ y_1, y_2 \} $ in $V$ such that the 
$\mathbb{R}$-linear map:
\begin{equation}
\label{orrr}
\langle \omega, \cdot\rangle  \colon V \otimes \mathbb{R} \rightarrow \mathbb{C} 
\end{equation}
is orientation preserving. We then define:
\begin{equation}
\label{deftau}
\tau(X,V) \ : = \ \frac{\langle \omega, y_1 \rangle}{\langle\omega, y_2 \rangle }. 
\end{equation} 
The quantity $ \tau(X,V)$ is a well-defined complex number which does not change
when one scales the class $\omega$. Due to the orientation preserving assumption
in $(\ref{orrr})$, its imaginary part $\tau_2(X,V)$ is positive. Moreover, 
if one modifies the oriented basis $ \{y_1, y_2\} $ on $V$, $ \tau(X,V)$ varies under the 
usual ${\rm SL}(2,\Zee) $-action on the upper half-plane $\mathbb{H}$. 
\par The second parameter mentioned earlier is introduced in the following form:
\begin{equation}
\label{defu2}
\widetilde{u}_2(X,V) \ : = \ 
\frac{\langle \omega, \overline{\omega} \rangle }{\vert \langle \omega, y_2 \rangle \vert^2 \cdot \ {\rm Im} \left [ \tau(X,V) \right ]}.
\end{equation}
\begin{lem}
\label{llleeemmm}
The term $ \widetilde{u}_2(X,V)$ is a positive real number and is independent of the choice of $\omega$. 
Moreover $ \widetilde{u}_2(X,V)$ remains 
unchanged when one modifies the oriented basis $ \{ y_1, y_2 \}$.  
\end{lem}
\begin{proof}
The first assertion is straightforward. To check the second assertion, note that if 
one changes the oriented basis in $V$ by way of an ${\rm SL}(2, \Zee)$ matrix 
$$ 
\left ( 
\begin{array}{cc}
a & b \\
c & d 
\end{array}
\right ),  
$$  
then this change induces the following variations:
$$ 
{\rm Im} \left [ \tau(X,V) \right ] \ \mapsto \ 
\frac{{\rm Im} \left [ \tau(X,V) \right ]}{\vert c \tau(X,V)+d\vert^2}
$$ 
$$ 
\vert \langle \omega, y_2 \rangle \vert ^2  \ \mapsto \ 
\vert c \tau(X,V)+d\vert^2 \cdot \vert \langle \omega, y_2 \rangle \vert ^2.
$$
The denominator of the right side of $(\ref{defu2})$ remains invariant. 
\end{proof}

\vspace{.05in}
\noindent Finally, we introduce the following auxiliary function:
$$ \rho \colon \mathbb{H} \rightarrow (0, \infty) $$
$$ \rho(\tau) \ = \ {\rm sup}_{m \in {\rm SL}(2, \Zee)} \ {\rm Im}\left ( m \cdot \tau \right ). $$
The function $ \rho(\tau)$ is well-defined and clearly invariant under the $ {\rm SL}(2, \Zee)$ action.
An alternative way of defining $ \rho(\tau) $ is as the imaginary part of the unique representative of 
the ${\rm SL}(2, \Zee)$ orbit of $\tau$ in the standard fundamental domain of $ \mathbb{H}$ with 
respect to the ${\rm SL}(2, \Zee) $ action. This formulation shows, in particular, that 
$ \rho(\tau) \geq \sqrt{3}/2 $.

\vspace{.07in}
\begin{dfn}
\label{vdef}
An elliptic $K3$ surface with section $ (X,\varphi,S) $ is said to have 
{\bf large complex structure with respect 
to} ${\bf V}$ if:
\begin{equation}
\label{llsscc1}
\widetilde{u}_2(X,V) \ > \ {\rm max} \left ( \ \rho (\tau(X,V)), \ \frac{2}{\sqrt{3}} \ \right ).
\end{equation} 
\end{dfn} 

\vspace{.05in}
\noindent In order to connect this definition with the discussion in section $\ref{rev}$,  
note that $\tau(X,V)$ and $\w{u}_2(X,V) $ can be 
seen as $C^{\infty}$ functions on $\Omega$ ($\tau$ is in fact analytic). One defines then 
the open subset $ \lcsd \subset \Omega $ of period lines satisfying condition $(\ref{llsscc1})$. 
This set decomposes as a disjoint union:
\begin{equation}
\lcsd \ = \ \lcsd^+ \ \sqcup \ \lcsd^- 
\end{equation}
where $\lcsd^+ = \lcsd \cap \Omega^+$ and $\lcsd^-$ is the complex conjugate of $\lcsd^+$. 
We claim then that $\lcsd^+$ satisfies the features 
described in $(\ref{lcsreq})$. Namely: 

\vspace{.05in}
\begin{theorem}
The open subset $\lcsd^+$ satisfies the following:
\begin{itemize}
\label{mainth}
\item [(a)] $\lcsd^+$ is left invariant by the action of $P^+$,
\item [(b)] $ {\rm U}({\rm N})_{\Zee} \backslash \lcsd^+$ is an open subset of the total space 
of the holomorphic principal $ \mathbb{T} $-bundle:
$$ {\rm U}({\rm N})_{\Zee} \backslash \Omega^+(V) \ \rightarrow \ {\rm U}({\rm N})_{\Cee} 
\backslash \Omega^+(V)  $$
defined in $(\ref{proj22})$ and its intersection with every fiber of the bundle is an open neighborhood of 
the cusp,
\item [(c)] for any two $ [\omega_1], [\omega_2] \in \mathcal{V}^+ $ and $\gamma \in \Gamma^+$ such that 
$ \gamma([\omega_1]) = [\omega_2] $, one has that $ \gamma \in {\rm P}^+ $.  
\end{itemize}
\end{theorem} 

\vspace{.03in}
\noindent The proof of Theorem $\ref{mainth}$ is presented in the next section.
\par We finish this section by formulating an extension of the large complex structure definition 
which is independent of ${\rm V}$. Recall from section $\ref{rev}$ that the set of primitive and 
isotropic rank-two sublattices in $ \left ( \mathcal{H}_X \right )^{\perp} $ is acted upon by 
the discrete group $\Gamma $ of isometries and that this action has two distinct orbits, essentially 
related to the two possible rank $16$ negative definite even and unimodular lattices
$E_8 \oplus E_8 $ and $ D_{16}^+$. As a consequence of Theorem $\ref{mainth} $, one has that, 
given an elliptic 
$K3$ surface with section $(X, \varphi,S)$ and two distinct isotropic lattices $V$ and $V'$ which belong 
to the same 
$\Gamma$-orbit, $(X, \varphi,S)$ can be of large complex structure with respect to at most one 
of $V$ and $V'$. This allows us to formulate: 
\begin{dfn}
\label{lcsdef}
An elliptic $K3$ surface with section $ X $ is said to have {\bf large complex structure in the} 
${\bf E_8 \oplus E_8 }$ $({\bf or} \ {\bf D^+_{16}})$ {\bf sense} if there exists a primitive 
isotropic rank-two lattice $V$ in the appropriate $\Gamma$-orbit such that $X$ has large complex 
structure with respect to $V$.
\end{dfn}   

\section{Proof of Theorem $\ref{mainth}$}   
\label{proof}     
As a first step, in order to gain a better understanding of the partial compactification outlined in section 
$\ref{rev}$, we introduce special coordinates on $\Omega^+$. We are going to use two distinct, 
but closely related, sets of parameterizations. The first parametrization is holomorphic and 
represents the standard {\bf tube domain realization} of $ \Omega^+$ (see, for example, 
\cite{dolgachev}). The second coordinate system is a (non-holomorphic) perturbation of the 
tube domain coordinates and presents $\Omega^+$ 
as a Siegel domain of the third kind. We shall refer to the latter coordinates as the 
{\bf Narain parametrization}, as the description defined by them is related to a string theory construction 
of Narain \cite{narain}. %Although not holomorphic in all directions, the Narain coordinates 
%have a good behaviour with respect to the action of the group of integral 
%isometries of ${\rm L}$. A description of this action plays an important role in our proof.    
\subsection{Narain Coordinates}
\label{narainsect}
As before, we start with a fixed choice of 
a rank-two sublattice ${\rm V}$ of ${\rm L}$ which is primitive and isotropic. For any 
$\omega$ representing a class in $\Omega$, the homomorphism:
\begin{equation}
\label{map21}
\langle \omega, \cdot \rangle \ \colon  {\rm V} \otimes \Ree \ \rightarrow \ \Cee 
\end{equation} 
is an isomorphism of real vector spaces. This fact follows immediately from the following lemma. 
\begin{lem}
\label{lem21}
Let $ v \in {\rm L} $ be a non-zero element satisfying $\langle v,v \rangle =0$. Then 
$ \langle \omega,v \rangle \neq 0 $ for any $[\omega] \in \Omega $.
\end{lem}
\bpf
Assume that $ \langle \omega, v \rangle =0$ for some $[\omega] \in \Omega$. Denote by $Q$ the plane 
in ${\rm L} \otimes \Ree$ spanned by the real and imaginary parts of $ \omega $. Then $Q$ is positive
definite with respect to $\langle \cdot, \cdot \rangle$ and since $v$ is not zero, $v \notin Q$. Then
$v$ and $Q$ span a three-dimensional subspace of  ${\rm L} \otimes \Ree$ on which the pairing  
$\langle \cdot, \cdot \rangle$ is non-negative. This contradicts the fact that the signature of ${\rm L}$ is 
$(2,18)$.
\epf

\noindent The above observation provides another effective method of differentiating between the two connected components of $\Omega$. 
If an orientation is chosen on $ {\rm V} \otimes \Ree$, then the map $(\ref{map21})$ is either orientation 
preserving or orientation reversing, depending on the component in which $[\omega]$ lies. 
\par For the purpose of streamlining future computations we shall fix at this point an orientation 
on $ {\rm V} \otimes \Ree $ such that $(\ref{map21})$ is orientation reversing for any $\omega$ representing 
a class in $\Omega^+$. We make then the first step toward parameterizing $\Omega^+$ by selecting a set $\{ x_1,x_2,y_1,y_2\} $ of four 
linearly independent isotropic elements of $ {\rm L}$ such that:
\begin{itemize}
\item $\{y_1,y_2 \} $ forms an oriented basis in $V$
\item $\langle x_1,x_2 \rangle  \ = \ \langle x_1,y_2 \rangle  \ = \ \langle x_2, y_1 \rangle \ = \ 0$
\item $ \langle x_1,y_1 \rangle  \ = \  \langle x_2,y_2 \rangle \ = \ 1$.
\end{itemize} 
It is clear that collections with the above features do exist. Such a collection is nothing but an embedding of an 
orthogonal direct sum $ H \oplus H $ into ${\rm L}$ in a way such the image contains ${\rm V}$. 
\par Denote then 
by $\Lambda$ the 
orthogonal complement in ${\rm L}$ of the sublattice generated by $x_1, x_2, y_1$ and $y_2$.  
It follows that $\Lambda$ is of rank $16$ and is unimodular and negative definite. Moreover, 
this construction induces a decomposition of ${\rm L}$ as a direct sum:
\begin{equation}
\label{ind11}
{\rm L} \ = \ 
\left ( \Zee x_1 \oplus \Zee x_2 \right ) 
\ \oplus  \ 
\left ( \Zee y_1 \oplus \Zee y_2 \right )
\ \oplus \
\Lambda.
\end{equation}
This allows us to identify an element of ${\rm L}$ as:
\begin{equation}
\label{ident666}
(a_1,a_2)(b_1,b_2)(c)
\end{equation}
where $a_1,a_2,b_1,b_2$ are integers, $c\in \Lambda$ and the parameters $a_1,a_2,b_1,b_2,c$ are uniquely 
determined. In this setting, the pairing 
$ \langle \cdot, \cdot \rangle $ is recovered as:
\begin{equation}
\left \langle \ (a_1,a_2)(b_1,b_2)(c), \ (a_1',a_2')(b_1',b_2')(c') \ \right \rangle \ = \ 
a_1b_1'+a_2b_2'+b_1a_1'+b_2a_2'+(c,c')
\end{equation}
where $(\cdot,\cdot)$ is the negative-definite pairing on $\Lambda$.
\begin{rem}
Note that the projection on the third term in $(\ref{ind11})$ induces a natural 
isometry:
$$ {\rm V}^{\perp}/{\rm V} \ \simeq \ \Lambda. $$
\end{rem}
\noindent Let then $[\omega] $ be an element of $\Omega^+$. Using Lemma $ \ref{lem21}$, one can always pick a normalized 
representative $\omega$ such that $\langle \omega, y_2 \rangle = 1$. Then, under the identification   
$(\ref{ind11})$, one can write:
\begin{equation}
\label{form1}
\omega \ = \ (\tau,1)(u,w)(z)
\end{equation}
where $\tau, u,v \in \Cee^* $, $z \in \Lambda_{\Cee}$. The first Hodge-Riemann 
condition $\langle \omega, \omega \rangle =0 $ becomes:
$$ 2 \left ( \tau u + w \right )+ (z,z) \ = \ 0 $$
and therefore:
\begin{equation}
\label{form2}
w \ = \ -\tau u - \frac{(z,z)}{2}. 
\end{equation}
After substituting $w$ in $(\ref{form1})$,  the second Hodge -Riemann condition 
$\langle \omega, \bar{\omega} \rangle > 0$ can be written as:
\begin{equation}
\label{cond21}
2 \tau_2 u_2 + (z_2,z_2)>0
\end{equation}
where the subscript 2 indicates that one takes the imaginary part. But 
$[\omega] \in \Omega^+$ and, because of the convention assuming that $(\ref{map21})$ reverses 
orientations, we obtain that $\tau_2 > 0$. Inequality $(\ref{cond21})$ 
implies then:
\begin{equation}
\label{ineq222}
u_2  \ > \  - \frac{(z_2,z_2)}{2 \tau_2} \ > \ 0.
\end{equation} 
It follows that:
\begin{prop} (Tube Domain Realization) One has an analytic isomorphism:
\begin{equation}
\label{tubemap}
\Omega^+ \ \stackrel{\simeq}{\longrightarrow} \ 
\left \{ 
\ (\tau,u,z) \in \mathbb{H} \times \mathbb{H} \times \Lambda_{\Cee} \ 
\vert \ 2 \tau_2 u_2 + (z_2,z_2)>0 \  
\right \}
\end{equation}
where $\mathbb{H} $ denotes the complex upper half-plane and $ \Lambda_{\Cee} = \Lambda \otimes \Cee$. 
The inverse of $(\ref{tubemap})$ maps a coordinate triple $(\tau,u,z)$ to a period line $[\omega]$ with:
\begin{equation}
\label{holoparam}
 \omega \ = \ (\tau,1) \left (u,  -\tau u - \frac{(z,z)}{2} \right ) (z).
\end{equation}  
\end{prop}
\noindent The {\bf Narain parametrization} perturbs slightly the second tube domain coordinate above. 
Indeed, under 
the same assumptions as before, let:
\begin{equation}
\w{u} \ \colon = \ u + \frac{(z,z_2)}{2\tau_2}. 
\end{equation}
The inequality $(\ref{ineq222})$ assures us that $\w{u} \in \mathbb{H} $ and it follows that:
\begin{prop} (Narain Coordinates) 
\label{naraincoord123}
One has a $ C^{\infty}$ isomorphism:
\begin{equation}
\label{narainmap}
\Omega^+ \ \stackrel{\simeq}{\longrightarrow} \ \mathbb{H} \times \mathbb{H} \times \Lambda_{\Cee} 
\end{equation}
which associates to a period line $[\omega]$ the triple $(\tau,\w{u},z)$. The inverse of $(\ref{narainmap})$ 
maps a coordinate triple $(\tau,\w{u},z)$ to a period line $[\omega]$ with:
$$ \omega \ = \ (\tau, \ 1) \left (\w{u}- \frac{(z,z_2)}{2\tau_2}, \  
-\tau   \w{u} + \frac{\tau (z,z_2)}{2\tau_2}  
 - \frac{(z,z)}{2} \right ) (z).$$  
\end{prop}
\noindent Note at this point that, for both parameterizations 
described above, the projection on the first and third 
coordinates:
\begin{equation}
\Omega^+ \ \rightarrow \ \mathbb{H} \times \Lambda_{\Cee}, \ \ [\omega] \mapsto (\tau,z)
\end{equation}
is holomorphic. This map is nothing but the restriction to $\Omega^+$ of the projection:
\begin{equation}
\Omega^+({\rm V}) \ \rightarrow \ {\rm U}({\rm N})_{\Cee} \backslash \Omega^+({\rm V})  
\end{equation}
from section $\ref{rev}$. Indeed, the integral nilpotent endomorphism ${\rm N}$ used 
in the construction of section $\ref{rev}$ can be written here explicitly as:
\begin{equation}
{\rm N} \colon {\rm L} \rightarrow {\rm L}, \ \ \ {\rm N}(\gamma) = 
\langle \gamma, y_2 \rangle y_1 -  \langle \gamma, y_1 \rangle y_2  
\end{equation} 
and, using the coordinates of $(\ref{ident666})$, it can be described as:
$$ {\rm N} \left ( (a_1,a_2)(b_1,b_2)(c) \right ) \ = \ (0,0)(a_2,-a_1)(0).$$
The group $U({\rm N})_{\Cee}$ acts then on a period line $[\omega]$ associated to 
$(\ref{holoparam})$ as follows:
\begin{equation}
{\rm exp}(w{\rm N}) \cdot [\omega] \ = \ \left [ 
(\tau,1) \left (u+w, -\tau (u + w) -  \frac{(z,z)}{2} \right ) (z)
\right ]. 
\end{equation}
Hence, in the framework of tube domain parametrization (as well as in Narain parametrization) the action 
of $U({\rm N})_{\Cee}$ has the effect of a translation in the second coordinate.
\par We also note that the two quantities $\tau(X,V)$ and $\w{u}_2(X,V)$ used in the definition of 
the large complex structure condition of section $\ref{largecond}$ are the first 
Narain coordinate and the imaginary part of the second Narain coordinate, respectively, of the corresponding 
period line in $\Omega^+$. The Narain coordinate description of the large complex structure domain is then:
$$ \lcsd^+ \ = \ \left \{ \ (\tau, \w{u}, z) \in \Omega^+ \ \vert \ \w{u}_2 > {\rm max} 
\left ( \rho(\tau), \frac{2}{\sqrt{3}} \right )   \right \}.$$

\subsection{The Action of the Parabolic Group ${\rm P}^+$}
Next, we use the framework constructed by the Narain coordinates to explicitly describe the action of 
the parabolic group ${\rm P}^+$ on the period domain $\Omega^+$. Note that, under 
decomposition $(\ref{ind11})$, one can view an isometry 
$ \gamma \in \Gamma $ as a matrix:
\begin{equation}
\label{izomatrix}
\left ( 
\begin{array}{ccc}
A & B & C \\
D & E & F \\
G & H & K \\
\end{array}
\right )
\end{equation}
with $ A,B,D,E $ in ${\rm End}(\Zee^2)$, $C,F$ in ${\rm Hom}(\Lambda, \Zee^2)$, $G,H$ in 
${\rm Hom}(\Zee^2, \Lambda)$ and $K \in {\rm O}(\Lambda)$. The conditions required for the 
entries of $ (\ref{izomatrix})$ to determine an actual isometry of ${\rm L}$ can be written 
explicitly and, they lead one to the following conclusion:
\begin{prop}
\label{teoparab}
A parabolic isometry $\gamma \in {\rm P}$ is given by a matrix:
\begin{equation}
\label{paramatrix}
\left ( 
\begin{array}{ccc}
m & 0 & 0 \\
R & \w{m} & - Qf \\
Q^tm & 0 & f \\
\end{array}
\right )
\end{equation}
with entries as follows.
\begin{enumerate}
\item $m \in {\rm GL}(2, \Zee)$ and $ \w{m} $ represents $ (m^t)^{-1}$.
\item $Q \in {\rm Hom}(\Lambda, \Zee^2)$, $R \in {\rm End}(\Zee^2)$ and they satisfy:
\begin{equation}
\label{rqcond}
R^tm + m^tR + m^tQQ^tm = 0. 
\end{equation}
\item $f$ is an isometry of $\Lambda$.
\end{enumerate}
A matrix as above corresponds to an isometry in ${\rm P}^+$ if and only if 
$m$ is an element of ${\rm SL}(2, \Zee)$. 
\end{prop}
\noindent The upper-script "t" refers to the adjoint of the homomorphism in question with 
respect to the pairings existing on its domain and target space. 
\par For simplicity, we shall refer to the isometry given by the matrix $(\ref{paramatrix})$ as 
$\gamma(m,Q,R,f)$. The composition law on ${\rm P}$ can then be read in this context as:
$$ \gamma(m_1,Q_1,R_1,f_1) \circ \gamma(m_2,Q_2,R_2,f_2) \ = \ 
\gamma(m_1m_2,\ Q_1+\w{m}_1Q_2f_1^{-1}, \ R_1m_2+\w{m}_1R_2+Q_1f_1Q_2^tm_2, \ f_1f_2).$$ 
In particular, we have:
\begin{equation}
\label{decomposit}
\gamma(m,Q,R,f) \ = \gamma(I,Q,m^{-1}R,I) \ \circ \ \gamma(m,0,0,I) \ \circ \ \gamma(I,0,0,f).
\end{equation}
\noindent This shows that ${\rm P^+}$ is generated by three special subgroups with familiar-looking structure. 
\begin{itemize}
\item [(a)] $ \mathcal{S} = \{ R=0, Q=0, f={\rm id}_{\Lambda} \} $. This group is naturally 
isomorphic to a copy of ${\rm SL}(2, \Zee)$.   
\item [(b)] $ \mathcal{W} = \{ m={\rm I}_2, R=0, Q= 0 \}$. This group is just the orthogonal group 
of the lattice $\Lambda$.
%\item [(c)] $ {\rm U}({\rm N})_{\Zee} = \{ m={\rm I}_2, Q= 0, f={\rm id}_{\Lambda} \} $. 
\item [(c)] $ \mathcal{T} = \{ m={\rm I}_2, f={\rm id}_{\Lambda} \} $. This is, essentially, the 
Heisenberg group of the lattice $\Lambda$. Note that, imposing $Q=0$ in $\mathcal{T}$, we obtain 
a normal subgroup which is no more but the abelian group $ {\rm U}({\rm N})_{\Zee} $. The quotient 
$ \mathcal{T}/{\rm U}({\rm N})_{\Zee}$ is naturally isomorphic to $ \Lambda \oplus \Lambda  $ and 
this fact leads to a presentation of $ \mathbb{T}$ as a semi-direct product:
$$ \left ( \Lambda \oplus \Lambda \right ) \ltimes 
{\rm U}({\rm N})_{\Zee}. $$  
\end{itemize}
By similar considerations, one verifies the following:
\begin{rem} The subgroups $\mathcal{S}$, $\mathcal{W}$ and $\mathcal{T}$ generate the entire ${\rm P}^+$. 
The following features also hold:
\begin{itemize}
\item [(1)] $ {\rm U}({\rm N})_{\Zee}  \subset {\rm Z} \left ( {\rm P}^+ \right ) $
\item [(2)] $ \mathcal{T} \vartriangleleft {\rm P}^+ $
\item [(3)] $ \left [ \mathcal{S} , \ \mathcal{W} \right ] \ = \ \{ I \} $ 
\item [(4)] $ {\rm P}^+ \ = \ \mathcal{T} \rtimes \left (\mathcal{S} \times \mathcal{W}  \right ) $.
\end{itemize}
\end{rem}
\noindent Let us then describe the action of ${\rm P}^+$ on $\Omega^+$ by writing the action for each 
of the three types of generators. We use Narain coordinates to make this description.

\vspace{.03in}
\begin{theorem} 
\label{parabtransf}
Let $[\omega]$ be a period line in $\Omega^+$ of Narain coordinates $(\tau,\w{u},z)$. Consider 
an integral isometry $\gamma \in {\rm P}^+$ and denote by $(\tau', \w{u}', z')$ the Narain 
coordinates of $ \gamma\left ( [\omega] \right )$. Then:
\begin{itemize}
\item [(a)] If $\gamma= \gamma(m,0,0,I) \in \mathcal{S}$ with $m$ given by the ${\rm SL}(2,\Zee)$
matrix:
$$ 
\left ( 
\begin{array}{cc}
a & b \\
c & d \\
\end{array}
\right )
$$ 
then:
\begin{equation}
\label{l1}
\tau' = \frac{a\tau+b}{c \tau + d}, \ \  \ \w{u}'=\w{u}, \ \ \ z' = \frac{z}{c\tau+d}. 
\end{equation}
\item [(b)] If $ \gamma= \gamma(I,0,0,f) \in \mathcal{W}$ with $f \in {\rm O}(\Lambda) $ 
then:
$$ \tau'=\tau, \ \ \ \w{u}' = \w{u}, \ \ \   z' = f(z). $$ 
\item [(c)] Assume that $\gamma = \gamma(I,R,Q,I) \in \mathcal{T} $. The homomorphism 
$Q \in {\rm Hom}(\Lambda, \Zee)$ induces two uniquely defined $c_1,c_2 \in \Lambda$ such 
that $ Q(c) = \left ( (c,c_1),(c,c_2) \right )  $ for any $c \in \Lambda$. Let 
$$ R \ = \ \left ( 
\begin{array}{cc}
r_{11} & r_{12} \\
r_{21} & r_{22} \\
\end{array}
\right ).$$
Then:
\begin{equation}
\label{l3}
\tau' =\tau, \ \ \ \w{u}' = \w{u} +  r_{12} + \frac{(c_1,c_1)}{2} - \frac{(z,c_1)}{2} + \frac{(\tau c_1 + c_2,z_2)}{2 \tau_2 }, 
\ \ \ z' = \tau c_1 + c_2 + z . 
\end{equation}
\end{itemize}  
\end{theorem}
\bpf
Let us analyze case (a). Note that: 
$$ \w{m} \ = \  \left ( 
\begin{array}{cc}
d & -c \\
-b & a \\
\end{array}
\right ). $$
Therefore, the parabolic transformation $\gamma(m,0,0,I)$ sends the period line:
$$ [\omega] \ = \ \left [ (\tau,1) \left (u, - \tau u -\frac{(z,z)}{2} \right ) (z) \ \right ] $$
to:
\begin{equation}
\label{sltransf}
\left [ (a\tau+b,c\tau+d) \left (du + c \tau u + \frac{c(z,z)}{2}, \ -bu -a \tau u - \frac{a(z,z)}{2} \right ) 
(z) \ \right ]. 
\end{equation}
After the proper normalization, the holomorphic (tube domain) coordinates of $ (\ref{sltransf}) $ 
are:
$$ \tau' = \frac{a\tau+b}{c \tau + d}, \ \  \ u'=u + \frac{c(z,z)}{2(c\tau+d)}, \ \ \ z' = \frac{z}{c\tau+d}. $$
We have then two of the conditions required by $(\ref{l1})$. The only thing missing is the second Narain coordinate of $ (\ref{sltransf}) $. But that value can be derived as:
\begin{equation}
\label{n1}
\w{u}' \ = \ u' + \frac{(z',z'_2)}{2\tau'_2} \ = \ u + \frac{c(z,z)}{2(c\tau+d)} + \frac{(z',z'_2)}{2\tau'_2} 
\ = \ \w{u} - \frac{(z,z_2)}{2 \tau_2} + \frac{c(z,z)}{2(c\tau+d)} + \frac{(z',z'_2)}{2\tau'_2}.
\end{equation}
The imaginary parts of the primed terms above are:
$$ \tau'_2 \ = \ \frac{\tau_2}{\vert c\tau+d\vert^2}, $$ 
$$ z'_2 \  = \ \frac{c\tau_1+d}{\vert c\tau+d\vert^2} \ z_2-\frac{c\tau_2}{\vert c\tau+d\vert^2} \ z_1  $$
and, after introducing these expressions in  $(\ref{n1})$ and carefully removing the canceling terms, we obtain that $\w{u}' \ = \  \w{u} $. Here the subscript 1 refers to the real part of the corresponding term.
\par Case (b) is straightforward. For case (c), we follow a strategy similar with the one used to check case (a).
The parabolic isometry $ \gamma(I,R,Q,I)$ sends the parabolic line: 
$$ [\omega] \ = \ \left [ (\tau,1) \left (u, - \tau u -\frac{(z,z)}{2} \right ) (z) \ \right ] $$
to
\begin{equation}
\label{qtransf}
\left [ (\tau,a) \left (r_{11} \tau + r_{12} + u - (z,c_1), \ r_{21}\tau + r_{22} -\tau u 
- \frac{(z,z)}{2} - (z,c_2) \right ) 
(z+ \tau c_1+ c_2 ) \ \right ]. 
\end{equation}
The tube domain coordinates of the above period line are then: 
$$ \tau' = \tau, \ \  \ u'=r_{11} \tau + r_{12} + u - (z,c_1), \ \ \ z' = z+ \tau c_1 + c_2 . $$
As with the case (a), what remains to be computed is the second Narain coordinate of $(\ref{qtransf})$. 
In order to evaluate this term, we write:
$$ \w{u}' \ = \ u' + \frac{(z',z'_2)}{2\tau'_2} \ = \  r_{11} \tau + r_{12} + u - (z,c_1) + \frac{(z+\tau c_1 +c_2 ,z_2 + \tau_2c_1)}{2\tau_2}   \ = \  $$
\begin{equation}
\label{zzzz}
\ = \ r_{11} \tau + r_{12} + \w{u} - \frac{(z,z_2)}{2\tau_2} - (z,c_1) + \frac{(z+\tau c_1 +c_2 ,z_2 + \tau_2c_1)}{2\tau_2} 
\end{equation}
Condition $(\ref{rqcond})$ imposes though that:
$$ r_{11} \ = \ - \frac{(c_1,c_1)}{2}. $$
After replacing $r_{11}$ in relation $(\ref{zzzz})$ and making the appropriate cancellation, one obtains:
$$ \w{u}'\ = \ \w{u} +  r_{12} + \frac{(c_1,c_1)}{2} - \frac{(z,c_1)}{2} + \frac{(\tau c_1 + c_2,z_2)}{2 \tau_2 }$$  
which is exactly the middle equality of $(\ref{l3})$.
\epf

\noindent Let us point out the following feature of the above result. Under parabolic transformations as in cases 
(a) and (b), the Narain coordinate $ \w{u}$ does not change. In case (c), $\w{u} $ is not invariant anymore.
However, its imaginary part $\w{u}_2$ is. Indeed, it follows from the middle equality of $(\ref{l3})$ that:
$$ \w{u}'_2 \ = \ \w{u}_2 - \frac{(z_2,c_1)}{2} + \frac{(\tau_2c_1,z_2)}{2 \tau _2} \ = \ \w{u}_2.$$
Since the three types of parabolic isometries of cases (a),(b), and (c) generate the entire group ${\rm P}^+$, we 
obtain that:    
\begin{cor}
The imaginary part $\w{u}_2$ of the second Narain coordinate is left invariant by the action of ${\rm P}^+$. 
\end{cor}
\noindent We note that the above corollary could also be derived from Lemma $\ref{llleeemmm}$.

\subsection{Proof of (a) and (b)}
From Theorem $\ref{parabtransf}$, we see that a parabolic transformation $\gamma \in {\rm P}^+ $,
acting on $\Gamma^+$, induces a ${\rm SL}(2,\Zee)$ change in the first Narain coordinate $\tau$ 
while preserving the imaginary part of the second Narain coordinate $\w{u}$. It follows that the 
quantities $\rho(\tau)$ and $\w{u}_2$ do not change under ${\rm P}^+$. The open subset 
$$ \lcsd^+ \ = \ \{ \ (\tau, \w{u}, z) \in \Omega^+ \ \vert \ \w{u}_2 > {\rm max} \left ( \rho(\tau), \frac{2}{\sqrt{3}} \right )   \}.$$
is therefore invariant under the action of ${\rm P}^+$.  
\par In order to verify (b), recall that the restriction to $\Omega^+$ of the holomorphic fibration:
$$
\Omega^+(V) \ \stackrel{p}{\rightarrow} \ {\rm U}({\rm N})_{\Zee} \backslash \Omega^+(V)
$$
is the projection in the first and last Narain coordinates:
$$ \Omega^+ \ \rightarrow \ \mathbb{H} \times \Lambda_{\Cee}, \ \ \ (\tau, \w{u}, z) \mapsto (\tau,z). $$
Since the action of the group $ {\rm U}({\rm N})_{\Zee} $ has as effect translations of $\w{u}$ by 
integers, we see that, under the exponential identification, $\lcsd^+ \cap p^{-1}(\tau,z)$ represents 
a punctured disc:
$$ \Delta^*_{(\tau,z)} \ = \ \{ \ t \in \Cee^* \ \vert \ \vert t \vert < e^{-2 \pi \cdot \ {\rm max}(\rho(\tau),
 2/\sqrt{3})} \ \}. $$
The center of the disc is the compactifying cusp.    
\subsection{Proof of (c)} 
Let $[\omega]$ and $ [\omega']$ be period lines in $\lcsd^+$. Assume that $\gamma \in \Gamma^+ $ has the 
the property that $ \gamma([\omega]) = [\omega']$. We shall prove here that, under these conditions, 
the isometry $\gamma$ must belong to ${\rm P}^+$.  
\par Most of the considerations required to check this fact will rely on the following technical argument:
\begin{lem}
\label{lemma1}
Let $r$ be an element of ${\rm L}$ which does not belong to ${\rm V}^{\perp}$. 
Denote by $a_1$ and $a_2$ the integers 
representing the intersection numbers $\langle r, y_1 \rangle$ and $\langle r, y_2 \rangle$, respectively. 
Then, for any $ [\omega] \in \Omega^+ $ with $\langle \omega, y_2 \rangle = 1$, one has: 
$$ \left \vert \left < \omega, r \right > \right \vert \ - \ 
\frac{\langle r,r \rangle \tau_2}{2\vert a_2 \tau - a_1 \vert}
\  \geq \ \vert a_2 \tau - a_1 \vert \cdot \w{u}_2  $$
where $\w{u}_2$ is the imaginary part of the second Narain coordinate of $[\omega]$. 
\end{lem} 
\begin{proof}
Let $(a_1,a_2)(b_1,b_2)(c)$ be the representation of $r$ in the decomposition $(\ref{ident666})$. 
Then, $a_1, a_2$ are not simultaneously zero and $2a_1b_1+2a_2b_2+(c,c)=\langle r,r \rangle $. 
\par The period line $[\omega]$ belongs to $\Omega^+$ and hence, it is given by a set of holomorphic 
(tube domain) coordinates $(\tau, u,z)$. For the sake of simplicity in the following computations, 
we shall denote $\langle \omega, r \rangle $ by $A$. We have:
$$ A \ = \ b_1\tau + b_2 + a_1u - a_2 \tau u - \frac{a_2}{2} (z,z) + (c,z)$$
and therefore,
\begin{equation}
\label{rel1}
\displaystyle{
u \ = \ \frac{b_1 \tau + b_2 - \frac{a_2}{2} (z,z) + (c,z)-A}{a_2 \tau - a_1}.
}
\end{equation}
The Narain coordinate $\tilde{u}$ is given by:
$$ \tilde{u} = u + \frac{(z,z_2)}{2 \tau_2}$$
and hence its imaginary part can be obtained as:
\begin{equation}
\tilde{u}_2 = u_2 + \frac{(z_2,z_2)}{2 \tau_2} = \frac{2\tau_2 u_2 + (z_2, z_2)}{2\tau_2}.
\end{equation}
But, from $(\ref{rel1})$, one computes:
\begin{equation}
\label{rel3}
 u_2 \ = \ \frac{{\rm Im}\left \{
\left [ b_1 \tau + b_2 - \frac{a_2}{2} (z,z) + (c,z) - A \right ] 
\left [ \overline{a_2 \tau - a_1} \right ] 
\right \} }{\vert a_2 \tau - a_1 \vert ^2} .  
\end{equation}
The numerator of the right-hand side of $(\ref{rel3})$ has the following form:
$$
{\rm Num} \ = \ \left [ b_1 \tau_1 + b_2 - \frac{a_2}{2} (z_1,z_1) + \frac{a_2}{2} (z_2,z_2) + (c,z_1) - A_1 \right ] 
\left [ -a_2 \tau_2 \right ] \ + $$
$$ +  \ 
\left [ b_1 \tau_2 - a_2 (z_1,z_2) + (c,z_2) - A_2 \right ] 
\left [ a_2 \tau_1 - a_1\right ] \ = 
$$
$$ = 
-[b_1a_1+b_2a_2]\tau_2 + \frac{a_2^2\tau_2}{2} (z_1,z_1) - \frac{a_2^2\tau_2}{2} (z_2,z_2)  \ - $$
$$ - \ 
a_2[a_2\tau_1-a_1](z_1,z_2) - a_2 \tau_2 (c,z_1) +  [a_2 \tau_1-a_1] (c, z_2) + 
$$ 
$$
+ A_1a_2\tau_2-A_2a_2\tau_1+ A_2a_1.
$$
Then:
\begin{equation}
\label{rel55}
\tilde{u}_2 \ = \ \frac{2\tau_2({\rm Num}) + \vert a_2 \tau - a_1 \vert ^2 (z_2,z_2) }{2\tau_2\vert a_2 \tau - a_1 \vert ^2}
\end{equation}
and one obtains the numerator of the right side fraction of $(\ref{rel55})$ as:
$$ 
-2[b_1a_1+b_2a_2]\tau_2^2 + a_2^2\tau_2^2 (z_1,z_1) - a_2^2\tau_2^2 (z_2,z_2)  \ + $$
$$  
- \ 2a_2[a_2\tau_1-a_1]\tau_2(z_1,z_2) - 2a_2 \tau_2^2 (c,z_1) +  2[a_2 \tau_1-a_1]\tau_2 (c, z_2)  \ + 
$$
$$ 
+ \ [a_2\tau_1-a_1]^2(z_2,z_2) +  a_2^2\tau_2^2 (z_2,z_2) \ +
$$
$$ 
+ \ 2A_1a_2\tau_2^2 - 2 a_2\tau_1\tau_2 +2 A_2 a_1 \tau_2 \ = 
$$
$$
= \ 
-2[b_1a_1+b_2a_2]\tau_2^2 + a_2^2\tau_2^2 (z_1,z_1) + [a_2\tau_1-a_1]^2(z_2,z_2) \ + \ $$
$$ - \ 
2a_2[a_2\tau_1-a_1]\tau_2(z_1,z_2) - 2a_2 \tau_2^2 (c,z_1) +  2[a_2 \tau_1-a_1]\tau_2 (c, z_2)  \  = 
$$
\begin{equation}
\label{rel4} 
-[2b_1a_1+2b_2a_2+(c,c)]\tau_2^2 \ + \ 
\end{equation}
$$ + \  \left ( [-a_2\tau_2]z_1+[a_2\tau_1-a_1]z_2+\tau_2c,  \ -[a_2\tau_2]z_1+[a_2\tau_1-a_1]z_2+\tau_2c 
\right ) \ + 
$$
$$ 
+ \ 2A_1a_2\tau_2^2 - 2 A_2a_2\tau_1\tau_2 +2 A_2 a_1 \tau_2. \
$$
Denote 
$$ w \colon = [-a_2\tau_2]z_1+[a_2\tau_1-a_1]z_2+\tau_2c. $$
We have then:
\begin{equation}
\tilde{u}_2 \ = \frac{
- \langle r,r \rangle \tau_2^2 + (w,w) + 2A_1a_2\tau_2^2 - 2 A_2a_2\tau_1\tau_2 +2 A_2 a_1 \tau_2
}{2\tau_2\vert a_2 \tau - a_1 \vert ^2}\ .
\end{equation}
One obtains therefore:
$$
\tilde{u}_2 - \frac{(w,w)}{2\tau_2\vert a_2 \tau - a_1 \vert ^2} \ = 
\frac{-\langle r,r \rangle \tau_2}{2\vert a_2 \tau - a_1 \vert ^2} \ + \ 
\frac{
2A_1a_2\tau_2^2 - 2 A_2a_2\tau_1\tau_2 +2 A_2 a_1 \tau_2
}{2\tau_2\vert a_2 \tau - a_1 \vert ^2} \ = 
$$
$$
 = \ \frac{-\langle r,r \rangle \tau_2}{2\vert a_2 \tau - a_1 \vert ^2} \ + \ \frac{
A_1a_2\tau_2 -  A_2a_2\tau_1 + A_2 a_1 
}{\vert a_2 \tau - a_1 \vert ^2} \ = 
$$
$$
= \  \frac{-\langle r,r \rangle \tau_2}{2\vert a_2 \tau - a_1 \vert ^2} \ + \  \frac{
{\rm Im} \left [ A(\overline{a_2\tau-a_1}  ) \right ] 
}{\vert a_2 \tau - a_1 \vert ^2} \ .
$$
Since the pairing $(\cdot, \cdot)$ is negative definite, it follows that:
\begin{equation}
\label{additionr}
\tilde{u}_2  \ \leq \  \frac{-\langle r,r \rangle \tau_2}{2\vert a_2 \tau - a_1 \vert ^2} \ + \  \frac{
{\rm Im} \left [ A(\overline{a_2\tau-a_1}  ) \right ] 
}{\vert a_2 \tau - a_1 \vert ^2} \ \leq \ \frac{-\langle r,r \rangle \tau_2}{2\vert a_2 \tau - a_1 \vert ^2} \ + \ 
\frac{\vert A \vert }{\vert a_2 \tau - a_1 \vert}.
\end{equation}
Multiplying the above line by $ \vert a_2 \tau - a_1 \vert$ produces the inequality stated in Lemma $\ref{lemma1}$.  
\end{proof}
\vspace{.2in} 
\begin{rem}
\label{remark11}
Let us remark that, during the course of the above proof, we actually show a slightly stronger 
result than the one stated in Lemma $\ref{lemma1}$. Namely, from the line just before $(\ref{additionr})$, it
follows that:
$$ \left \vert {\rm Im} \left ( \frac{\langle \omega,r \rangle }{a_2\tau-a_1 }\right ) \right \vert \ -
 \ \frac{\langle r,r \rangle \tau_2}{2\vert a_2 \tau - a_1 \vert ^2} \ \geq \ \widetilde{u}_2. $$
\end{rem} 
\vspace{.2in} 
\noindent Let us return to the proof of $\ref{mainth}$ (c). Our strategy is as follows. In holomorphic (tube domain) coordinates: 
$$ \omega \ = \ (\tau,1) \left (u, -\tau u - \frac{1}{2}(z,z) \right )(z) $$ 
$$ \omega' \ = \ (\tau',1) \left (u', -\tau' u' - \frac{1}{2}(z',z') \right )(z'). $$
It would suffice to prove $\ref{mainth}$ (c) for the case when 
\begin{equation}
\label{partcond}
\rho(\tau) = \tau_2 \ \ {\rm and} \ \  \rho(\tau')= \tau_2', 
\end{equation}
as, from the description in Theorem $\ref{parabtransf}$, 
one can always find parabolic isometries in $\mathcal{S} $ that transform $\omega$ and $\omega'$ to 
period lines satisfying the above conditions.  
\par We shall therefore assume that $(\ref{partcond})$ holds. We prove then that, under these conditions, 
if $\gamma([\omega]) = [\omega']$, then $ \gamma^{-1}(y_1) $ and $ \gamma^{-1}(y_2) $ must belong to $V$. This claim implies $\gamma^{-1} \in {\rm P}^+$ which in turn gives $\gamma \in {\rm P}^+$. 

\vspace{.03in}
\begin{claim}
\label{lemma2}
$\gamma^{-1}(y_2) \in V$. 
\end{claim}
\bpf
\noindent Since $ \gamma([\omega]) = [\omega']$, one can write that, for some $ \alpha \in \Cee^*$,
\begin{equation}
\label{rel77}
\alpha \left  ( \ (\tau',1)(u', -\tau u' - \frac{1}{2}(z',z'))(z') \ \right ) \ = \ 
\gamma \left ( \ (\tau,1)(u, -\tau u - \frac{1}{2}(z,z))(z) \ \right ). 
\end{equation}
The scaling factor $ \alpha $ can be determined as follows:
$$ \alpha \ = \ \langle \gamma(\omega), y_2 \rangle \ = \ \langle \omega, \gamma^{-1}(y_2) \rangle. $$
Since $\gamma$ is an isometry of ${\rm L}$, $\gamma^{-1}(y_2)$ is integral primitive and isotropic. 
In the framework of $(\ref{ident666})$, one can then write:  
$$ \gamma^{-1}(y_2) = (a_1,a_2)(b_1,b_2)(c) $$
with $ 2(a_1b_1+a_2b_2)+(c,c)=0$. 
\par Recall then that $\langle \omega, \overline{\omega} \rangle \ = \ 2 \tau_2 u_2 + (z_2, z_2) \ = \ 2 \tau_2 \tilde{u}_2$. This feature, in connection with $ (\ref{rel77})$ provides:
\begin{equation}
\tau_2 \tilde{u}_2 \ = \ \frac{\langle  \omega, \overline{\omega} \rangle }{2}  \ = \frac{\vert \alpha \vert^2 
\langle \omega', \overline{\omega}' \rangle }{2} \ = \ \vert \alpha \vert ^2 \tau'_2 \tilde{u}'_2.   
\end{equation}  
Hence:
\begin{equation}
\label{rel88}
\tau'_2 \tilde{u}'_2 \ = \ \frac{\tau_2 \tilde{u}_2}{\vert \alpha \vert ^2}.
\end{equation}
We claim that the above condition implies that $a_1$ and $a_2$ are simultaneously zero. 
\par Let us assume the contrary and show that a contradiction follows. 
Lemma $\ref{lemma1}$ applied for $r=\gamma^{-1}(y_2)$ implies that:
$$ \vert \alpha \vert  \ \geq \ \vert a_2\tau - a_1 \vert \tilde{u}_2.  $$
Therefore, 
\begin{equation}
\label{cond9999}
\tau'_2 \tilde{u}'_2 \ \leq \ \frac{\tau_2}{\vert a_2\tau- a_1 \vert ^2\tilde{u}_2 }.
\end{equation}
But, since the assumption is that $a_1$ and $a_2$ are not both zero, one has that:
$$ \frac{\tau_2}{\vert a_2\tau- a_1 \vert ^2  } \ \leq \ \rho (\tau). $$ 
Then $(\ref{cond9999})$ and the fact that $[\omega] \in \lcsd^+$ lead to:
$$ \tau'_2 \tilde{u}'_2 \  \leq \
\frac{\rho(\tau)}{\tilde{u}_2} \ < \ 1. $$
One has then that:
$$ \tilde{u}'_2 \ < \ \frac{1}{\tau'_2} = \frac{1}{\rho(\tau')} \ \leq \ \frac{2}{\sqrt{3}} $$
which clearly contradicts $ [\omega'] \in \lcsd^+$.
\par We have therefore that $ a_1=a_2=0$. Then, since $\gamma^{-1}(y_2)$ is isotropic, 
we have that $(c,c)=0$ which, in turn, implies that $c=0$. It follows that $ \gamma^{-1}(y_2)\in V$. 
\epf  \

\vspace{.03in}
\noindent Next, we complete the final step in the proof of $\ref{mainth}$(c).
\vspace{.02in}
\begin{claim}
\label{lemma3}
$\gamma^{-1}(y_1) \in V$.
\end{claim}
\bpf 
Claim $\ref{lemma2}$ assures us that: 
$$ \gamma^{-1}(y_2) = (0,0)(b_1,b_2)(0)$$
where $b_1,b_2$ are not simultaneously zero and are relatively prime. Let us then analyze the integral 
element:
$$ \gamma^{-1}(y_1) \ = \ (a'_1,a'_2)(b'_1,b'_2)(c').$$
If $a'_1=a'_2=0$, this fact together with the fact that $\gamma^{-1}(y_1)$ is isotropic implies 
that $c'=0$. It follows then that $ \gamma^{-1}(y_1) \in {\rm V}$ and the proof of Claim $ \ref{lemma3}$ is done. 
\par It suffices therefore to check that the assumption that $a'_1 $ and $ a'_2$ cannot hold. Indeed, we show that 
by assuming such a case, one is led to a contradiction. Since:
$$ \langle \gamma^{-1}(y_1), \gamma^{-1}(y_2) \rangle \ = \ \langle y_1,y_2 \rangle \ = \ 0 $$
we have that :
\begin{equation}
\label{eq11}
a'_1b_1+a'_2b_2=0.
\end{equation}
Relation $(\ref{eq11})$ combined with the fact that $b_1, b_2$ are relatively prime implies then 
that:
$$ a'_1 = q b_2, \ a'_2 = -qb_1$$
for some $q \in \Zee$, $q \neq 0$.
\par But then, recalling equation $(\ref{rel77})$ of Claim $\ref{lemma2}$, we deduce that:
$$ \tau' \ = \ \frac{\langle \omega, \gamma^{-1}(y_1) \rangle }{\langle \omega, \gamma^{-1}(y_2) \rangle } \ = \ 
\frac{\langle \omega, \gamma^{-1}(y_1) \rangle }{b_1\tau+b_2}. $$
This implies:
\begin{equation}
\label{cond2323}
\tau' \ = \ q \frac{ \langle \omega, \gamma^{-1}(y_1) \rangle }{a'_2\tau-a'_1}. 
\end{equation}
At this point, Remark $ \ref{remark11} $ applied to $ r = \gamma^{-1}(y_1) $ and $\omega$ 
provides the following estimate:
$$ 
\left \vert {\rm Im} 
\left ( 
\frac{\langle \omega, \gamma^{-1}(y_1)\rangle}{a'_2\tau-a'_1} 
\right )
\right \vert  \ \geq \ \w{u}_2. $$ 
From $(\ref{cond2323})$, one obtains then: 
\begin{equation}
\label{cond24242}
 \tau'_2 \ = \ \vert q \vert \cdot \left \vert {\rm Im} 
\left ( 
\frac{ \langle \omega, \gamma^{-1}(y_1) \rangle }{a'_2\tau-a'_1} 
\right )
\right \vert \ \geq \ \vert q \vert \cdot \widetilde{u}_2 .
\end{equation}  
But, equation $(\ref{rel88})$ tells us that:
\begin{equation}
\label{eq33}
\tau'_2\widetilde{u}'_2 \ = \ \frac{\tau_2\widetilde{u}_2}{\vert b_1 \tau + b_2 \vert ^2 }.
\end{equation}
This fact, together with inequality $ (\ref{cond24242})$ provides:
\begin{equation}
\label{eq22}
\frac{\tau_2}{\vert b_1\tau_2+b_1\vert^2 \widetilde{u}'_2 } \ \geq \ \vert q \vert. 
\end{equation}
Let us analyze the possibilities that can appear here. If one assumes: 
$$ 
\frac{\tau_2}{\vert b_1\tau+b_2\vert^2} \ \leq \ 1
$$
then inequality $ (\ref{eq22})$ implies that 
$$ \vert q \vert \ \leq \ \frac{1}{\widetilde{u}'_2} \ < \ \frac{1}{2} $$
which contradicts the fact that $q$ is integral and not zero.
Therefore, it must be the case that:
$$ 
\frac{\tau_2}{\vert b_1\tau+b_2\vert^2} \ > \ 1.
$$
However, since $\tau_2 = \rho(\tau)$, $\tau_2$ represents the maximum imaginary part over 
the ${\rm SL}(2, \Zee)$ orbit of $\tau$. The above condition can then only hold if $b_1=0$ and 
$b_2=1$. In such a situation, taking into account $ (\ref{eq22})$ and $ (\ref{eq33})$, one writes:
$$ 1 \ \leq \ \frac{\tau'_2}{\widetilde{u}_2} \ = \ \frac{\tau_2}{\widetilde{u}'_2} .$$
Therefore:
$$ \tau'_2 \ \geq \w{u}_2 \ \ {\rm and} \ \ \tau_2 \geq \widetilde{u}'_2. $$
But the above inequalities contradict the fact that $[\omega]$ and $[\omega']$ are in 
$ \lcsd^+$, which implies:
$$ \widetilde{u}_2 > \tau_2 \ \ {\rm and} \ \ \widetilde{u}'_2 > \tau'_2. $$   
We conclude therefore that the assumption that $a'_1,a'_2$ are not simultaneously zero cannot hold. It 
follows that $a'_1=a'_2=0$, and, as explained earlier, this implies $\gamma^{-1}(y_1) \in {\rm V}$. 
\epf

\section{A Consequence of the Large Complex Structure Condition}
\label{ade}
In this section we analyze the types of singular fibers (or rather their so-called ADE types) that appear 
in the elliptic fibration of a triple $(X,\varphi,S)$ satisfying the large complex structure condition. 
Recall that, if 
$$ \mathcal{U} \ \subset \ \Gamma \backslash \Omega $$
is one of the two large complex structure domains defined in section $\ref{largecond}$, there exists
a natural holomorphic fibration with fibers isomorphic to open punctured discs:
\begin{equation}
\label{largefib}
\mathcal{U} \ \rightarrow \ \mathcal{D} 
\end{equation}
whose base space is the appropriate Type II Mumford boundary divisor. We shall prove here that the 
ADE type does not change over the fibers of $(\ref{largecond})$. 
\par In order to position the above statement on a rigorous footing, let us review a few 
classical definitions and results. For details and explicit 
proofs we refer the reader to \cite{kodaira1} \cite{bpv} and \cite{miranda2}. 
\par Let $(X, \varphi, S)$ be an elliptic K3 surface with section. As mentioned earlier, in such a context 
one has a decomposition of the Neron-Severi lattice:
$$ {\rm NS}(X) \ = \mathcal{H}_X \oplus \mathcal{W}_X $$
where $ \mathcal{W}_X $ is the negative-definite sublattice of ${\rm NS}(X)$ generated by classes 
associated to algebraic cycles orthogonal to both the elliptic fiber and the section. 

\vspace{.03in}
\begin{dfn}
The sublattice $\mathcal{W}_X^{{\rm root}}$ of $ \mathcal{W}_X$ spanned by:
$$ \{ \ r \in \mathcal{W}_X \ \vert \ \langle r,r \rangle = -2 \ \} $$
is called the {\bf ADE type} of the elliptic fibration with section $(X,\varphi,S)$.
\end{dfn}

\vspace{.03in}
\noindent The reason for the above terminology is that, the lattice $\mathcal{W}_X^{{\rm root}}$ has a special decomposition involving the classical root lattices $A_n$ $D_n$ and $E_n$ and, this decomposition encodes important information about the geometry of the singular fibers of the elliptic pencil $\varphi$. In order 
to explain this feature, denote by $\Sigma$ the finite set of points $v \in \mathbb{P}^1 $ for which the corresponding fiber $F_v$, in the elliptic fibration $ \varphi \colon X \rightarrow \mathbb{P}^1 $, is singular. 
For each $v \in \Sigma$ one has a formal decomposition into irreducible components:
\begin{equation}
\label{adedecompo}
F_v \ = \ \Theta_{v,0} \ + \ \sum_{j=1}^{t_v-1} \ \mu_{v,j} \Theta_{v,j}.
\end{equation}   
Here $t_v \geq 1$ represents the number of irreducible components of $F_v$ and $\Theta_{v,0}$ is the unique 
irreducible component of $F_v$ meeting the section $S$. One denotes then by $T_v$ the sublattice in  
$\mathcal{W}_X$ spanned by the classes:
$$ c_1 \left ( \Theta_{v,j} \right ), \ \ 1 \leq j \leq t_v-1. $$
The following classical result due to Kodaira \cite{kodaira1} relates the lattices $T_v$ with the 
geometry of the singular fibers.

\vspace{.03in} 
\begin{theorem} 
\label{kodairaclasif}
If $t_v \geq 2$ then $ \Theta_{v,j}$ is a smooth rational curve for $0\leq j \leq t_v-1$. Moreover, one can 
deduce the isomorphism class 
of the lattice $T_v$ from the Kodaira type of the singular fiber $F_v$ as follows:
$$ \begin{tabular}{|r|r|}
\hline
{\rm Type \ of} $F_v$& $T_v$ \\
\hline
$I_1$, $II$&\{0\}\\
$I_2$, III &$A_1$\\
$I_3$,  IV &$A_2$\\
$I_{n}$  (n$\geq$4)&$A_{n-1}$\\
$I^*_{n}$&$D_{n+4}$\\
$IV^*$&$E_6$\\
$III^*$&$E_7$\\
$II^*$&$E_8$\\
\hline
\end{tabular}
$$
\end{theorem} 

\vspace{.03in}
\noindent Note at this point that $T_v \subset \mathcal{W}^{{\rm root}}_X$. Also, for $v_1 \neq v_2$, the two 
lattices $T_{v_1}$ and $T_{v_2}$ are orthogonal. This allows 
one to define the direct sum:
\begin{equation}
\label{adedecomp}
T \ = \ \bigoplus _{v \in \Sigma} \ T_v \ \subset \ \mathcal{W}^{{\rm root}}_X.
\end{equation} 
\begin{prop}
$ T = \mathcal{W}^{{\rm root}}_X$.
\end{prop}
\bpf
It suffices to check that any root of $\mathcal{W}_X$ also belongs to $T$. Let $r$ be such a root. 
Since $\langle r, r \rangle = -2$, by the Riemann-Roch theorem, either $r$ or $-r$ represents 
an effective divisor $D$ on $X$. Let: 
$$ D = \sum n_i D_i $$ 
be the formal decomposition of $D$ into irreducible 
components and denote by $F$ the divisor class of the elliptic fiber. Since $F \cdot D=0$ and $F$ is nef, we 
deduce that $F$ has vanishing intersection with each of the irreducible components $D_i$. It follows 
that each $D_i$ is either equivalent to $F$, or it is an irreducible component of a singular fiber $F_v$ 
for some $v \in \Sigma $. But $S\cdot D=0$ and, a simple look at the decomposition 
$ (\ref{adedecompo})$ assures us that $D_i \in T_v$.    
\epf \

\noindent We have therefore a decomposition:
\begin{equation}
\label{adeedecomp}
\mathcal{W}^{{\rm root}}_X \ = \ \bigoplus _{v \in \Sigma} \ T_v
\end{equation}
in which every term is isomorphic to one of the classical root lattices $A_n$ ($n \geq 1$), $D_n$ ($n \geq 1$) 
or $E_n$ ($n=6,7,8$). Moreover, since the root lattices are known to be indecomposable, the decomposition 
$(\ref{adeedecomp})$ is unique. Therefore, by knowing the isomorphism class of $\mathcal{W}^{{\rm root}}_X$, one 
is able to detect, in a lattice-theoretic manner, most of the geometric types of singular fiber appearing 
in the actual elliptic fibration. Of course, this method does not distinguish between the Kodaira 
types $I_2$ and III or between $I_3$ and IV and also cannot detect the 
appearance of the singular fiber types  $I_1$ or $II$.
\par Let us formulate the result suggested at the beginning of the section. As in the construction 
of section $\ref{rev}$, we assume a choice of a rank-two primitive isotropic sublattice 
$ {\rm V} \subset {\rm L}$. The associated Type II Mumford boundary component is:
$$ \mathcal{D}(V) \ = \ {\rm P}^+ \backslash \left ( {\rm U}(N)_{\Cee} \backslash \Omega^+({\rm V}) \right ). $$  
The large complex structure domain associated to ${\rm V}$
$$ \mathcal{U} \ \subset \ \Gamma \backslash \Omega $$
fibers holomorphically:
\begin{equation}
\label{fibbb}
\mathcal{U} \ \rightarrow \ \mathcal{D}({\rm V})
\end{equation}
with all fibers being isomorphic to punctured complex planes.

\vspace{.04in}
\begin{theorem}
\label{conseq}
The ADE type lattice $\mathcal{W}^{{\rm root}}_X$ remains constant over the fibers of $(\ref{fibbb})$. 
\end{theorem}
\bpf
It suffices to prove the above statement for the projection:
\begin{equation}
\label{fibbb1}
\mathcal{V}^+ \ \rightarrow \  {\rm U}(N)_{\Cee} \backslash \Omega^+({\rm V}) 
\end{equation}
which covers $(\ref{fibbb})$. But, as we explained during the ending remarks of section $\ref{narainsect}$, 
in the framework of Narain coordinates, $(\ref{fibbb1})$ is just:
\begin{equation}
\label{narprojj}
 \{ \ (\tau, \w{u}, z) \ \vert \ \w{u}_2 > {\rm max}\left ( \rho(\tau),\frac{2}{\sqrt{3}} \right ) \ \} \ \rightarrow \ 
 \mathbb{H} \times \Lambda_{\Cee} 
\end{equation}   
$$ (\tau, \w{u},z) \ \rightarrow \ (\tau,z). $$ 
Let then
$$ \mathcal{W}^{{\rm root}}_{(\tau, \w{u},z)} \ \subset \ {\rm L} $$ 
be the ADE lattice associated to the period line $[\omega] \in \Omega^+$ of Narain 
coordinates $(\tau,\w{u},z)$. Theorem $\ref{conseq}$ is implied by the following two 
claims.
\begin{enumerate} 
\item  $  \mathcal{W}^{{\rm root}}_{(\tau, \w{u},z)} \ \cap \ {\rm V}^{\perp} \ $ does not depend on $\w{u}$.
\item If $\w{u}_2> {\rm max} \left ( \rho(\tau), \frac{2}{\sqrt{3}}  \right ) $, then $  \ \mathcal{W}^{{\rm root}}_{(\tau, \w{u},z)} \subset {\rm V}^{\perp} $. 
\end{enumerate} 
\noindent The first claim is almost straightforward. Recall the decomposition $(\ref{ident666})$ and the 
Narain parametrization of Proposition $\ref{naraincoord123}$. An element of ${\rm L}$  
$$ r \ = \ (a_1,a_2)(b_1,b_2)(c)  $$
is a root in  $ \mathcal{W}^{{\rm root}}_{(\tau, \w{u},z)} $ if and only if the following two conditions hold.
\begin{equation}
\label{cond1}
\langle r, r \rangle \ = \ 2(a_1b_1+a_2b_2)+(c,c)\ = \ -2
\end{equation}
\begin{equation}
\label{cond2}
\langle \omega, r \rangle \ = \ b_1\tau+b_2 +a_1 \left ( \w{u}-\frac{(z,z_2)}{2\tau_2} \right ) + a_2 
\left ( 
-\tau \w{u} + \frac{\tau(z,z_2)}{2\tau_2} - \frac{(z,z)}{2}
\right )   
+ (c,z) = 0.
\end{equation}
But ${\rm V}^{\perp}$ corresponds to $a_1=a_2=0$. Therefore:
$$  \mathcal{W}^{{\rm root}}_{(\tau, \w{u},z)} \ \cap \ {\rm V}^{\perp} \ = \ \left \{ \  
 (0,0)(b_1,b_2)(c) \ \vert \ (c,c)=-2, \ b_1\tau + b_2 + (c,z)=0 \
\right \}$$
and clearly it does not depend on $\w{u}$.
\par In order to justify the second claim, we show that, under the assumption 
$\w{u}_2>{\rm max} \left ( \rho(\tau), 2/\sqrt{3} \right ) $, 
conditions $(\ref{cond1})$ and $(\ref{cond2})$ imply that $a_1=a_2=0$. To check this, 
let us assume that $a_1$ and $a_2$ are not simultaneously zero. Then, if 
conditions $(\ref{cond1})$ and $(\ref{cond2})$ are satisfied, lemma $\ref{lemma1}$ tells us that:
\begin{equation}
\label{rel998}
\w{u}_2 \ \leq \ \frac{\tau_2}{\vert a_2 \tau - a_1 \vert^2}.
\end{equation}
\begin{itemize}
\item [(a)] If $a_1 \neq 0$ and $a_2 \neq 0$, let $n={\rm gcd}(\vert a_1 \vert ,\vert a_2 \vert )$. Inequality $(\ref{rel998})$ implies 
then 
$$ \w{u}_2 \ \leq \ \frac{\rho(\tau)}{n^2} \ \leq \ \rho(\tau). $$  
\item [(b)] If $a_1 = 0$ and $a_2 \neq 0$ then $(\ref{rel998})$ gives:
$$ \w{u}_2 \ \leq \ \frac{\tau_2}{\vert a_2 \vert ^2 \vert \tau \vert ^2 } \ \leq \ \frac{\rho(\tau)}{\vert a_2 \vert ^2} \ \leq \ \rho(\tau). $$ 
\item [(c)] If $a_1 \neq 0$ and $a_2 = 0$ then, from $(\ref{rel998})$:
$$ \w{u}_2 \ \leq \ \frac{\tau_2}{\vert a_1 \vert ^2  } \ \leq \ \frac{\rho(\tau)}{\vert a_1 \vert ^2} \ \leq \ \rho(\tau). $$ 
\end{itemize}
All three possible cases produce contradictions with the large complex structure assumption 
$\w{u}_2 > \rho(\tau)$.  
\epf \

\noindent We close this section with a final comment about the above result. We have shown that, if 
one moves inside the moduli space $ \Gamma \backslash \Omega $, from a Type II boundary point and 
following the fibers of $(\ref{fibbb})$, the ADE type of the elliptic pencils associated to the 
points encountered stays constant. The fibers of $(\ref{fibbb})$ are coming from nilpotent orbits 
on $\Omega^+$ and, therefore, they can be continued indefinitely inside $ \Gamma \backslash \Omega $. 
However, outside $\mathcal{U}$ they no longer form the fibers of a fibration over the boundary divisor. 
The images of nilpotent orbits corresponding to distinct boundary points may intersect inside the 
moduli space. 
\par On each nilpotent orbit, one encounters points where the lattice $\mathcal{W}^{{\rm root}}_X$ 
contains roots no longer belonging to ${\rm V}^{\perp}$. At such points, some exceptional 
singular fibers appear in the corresponding elliptic pencils, as the ADE type lattice is no longer 
generic (within the given nilpotent orbit). 
These are the points physicist refer to as points of {\bf enhanced gauge symmetry}. The 
reason for this terminology is that, under the eight-dimensional F-theory/heterotic string duality, these
points correspond on the heterotic side to flat $G$-bundles with reduced structure group.  
\par In light of this interpretation, Theorem $\ref{conseq} $ states that all enhanced gauge 
symmetry points lie outside of our large complex structure region.

\end{document}